\documentclass[twoside,12pt]{article}
\usepackage{amsmath}
\usepackage{amssymb}
\usepackage[usenames,dvipsnames]{color}
\def\takeshi{\color{red}}           
\def\pier{\color{blue}}               
\def\revis{\color{Green}}            
\def\colli{\color{blue}}
\def\fukao{\color{red}}
\let\takeshi\relax
\let\pier\relax
\let\revis\relax
\let\colli\relax
\let\fukao\relax

\title{The Allen-Cahn equation\\ with dynamic boundary conditions\\ 
and mass constraints
}
\author{Pierluigi Colli\\
Dipartimento di Matematica, Universit\`a degli Studi di Pavia\\
Via Ferrata~1, 27100 Pavia, Italy\\
E-mail: \texttt{pierluigi.colli@unipv.it}\\
\and \\ Takeshi Fukao\\
Department of Mathematics, Faculty of Education\\
Kyoto University of Education\\
1~Fujinomori, Fukakusa, Fushimi-ku, Kyoto~612-8522 Japan\\
E-mail: \texttt{fukao@kyokyo-u.ac.jp}}

\newcommand\testopari{\sc Pierluigi Colli and Takeshi Fukao}
\newcommand\testodispari{\sc Allen-Cahn equation with dynamic b.c.\ and mass constraints}
\begin{document}

\date{}

\maketitle

\begin{abstract}
{\pier The Allen-Cahn equation, coupled with dynamic boundary conditions,
has recently received a good deal of attention. The new issue of 
this paper is the setting of a rather general mass constraint 
which may involve either the solution inside the domain 
or its trace on the boundary. The system of nonlinear partial 
differential equations can be formulated as variational inequality.
The presence of the constraint in the evolution process leads to 
additional terms in the equation and the boundary condition 
containing a suitable Lagrange multiplier. A well-posedness result is 
proved for the related initial value problem.}

\vspace{2mm}

\vspace{2mm}
\noindent \textbf{AMS (MOS) subject clas\-si\-fi\-ca\-tion:} 35K86, 49J40, 80A22.

\end{abstract}

\section{Introduction}
\setcounter{equation}{0}

{\pier The Allen-Cahn equation \cite{AlCa} is a famous equation aiming to {\revis describe} the order-disorder  phase transition in a process of phase separation in a binary alloy. It is applicable to several directions, for example, it is widely employed in the description of the solid-liquid phase transition (see the monograph \cite{BrSp} and references therein).

Let $0<T<+\infty$ and $\Omega \subset \mathbb{R}^{\takeshi d}$, 
${\takeshi d}=2$ or $3$, be the bounded smooth domain 
occupied by the material. 
Also the boundary $\Gamma$ 
of $\Omega $ is supposed to be smooth enough. 
We recall the isothermal Allen-Cahn equation in  
the following form:
\begin{equation}
	\frac{\partial u}{\partial t}-\Delta u +W'(u) = f
	\quad \mbox{a.e.\ in }Q:=\Omega \times (0,T), 
\label{i1}
\end{equation}
where the unknown $u:=u(x,t)$ stands for the order parameter
and $f:=f(x,t)$ is a given source term.  
The nonlinear term $W'$ plays an important role, it is the derivative 
of a function $W$ {\revis usually} referred as double well potential,
with two minima and a local unstable maximum in between. 
The prototype model for the Allen-Cahn equation 
is provided by $W(r)=(1/4)(r^2-1)^2$, $r\in \mathbb{R}$,
so that $W'(r)=r^3-r$, $r\in \mathbb{R}$, is the sum of 
an increasing function with a power growth and another smooth 
(in particular, Lipschitz continuous)
function which breaks the {\revis monotonicity} properties of 
the former (and is related to the non-convex part 
of the potential $W$). In this paper, we treat more general cases
for such a nonlinearity, that is, we assume that $W'$ is the sum
of a maximal monotone graph (it can be a graph with vertical segments too)
defined in the whole of  $\mathbb{R}$ and of a Lipschitz perturbation. 

Usually, the Allen-Cahn equation is coupled with the homogeneous Neumann boundary 
condition, which means no flux exchange at the boundary.  
Recently, equation~\eqref{i1} has been investigated 
(see, e.g., {\takeshi \cite{CaCo, CoSp, GaGu, GaWa, Li}} and references therein) when complemented by a dynamic boundary condition of the 
following form:
\begin{equation}
	{\takeshi u_{\Gamma}=u_{|_\Gamma}} , \quad 
	\partial_\nu u 
	+ \frac{\partial u_{\Gamma}}{\partial t}-\kappa \Delta_{\Gamma } u_{\Gamma}
	+W_{\Gamma }'(u_{\Gamma}) = f_\Gamma \quad \mbox{a.e.\ on }\Sigma:= \Gamma \times (0,T).
\label{i2}
\end{equation}
Here, {\colli $u_{|_\Gamma}$ denotes the trace of $u$ and $\partial_\nu$ represents} the outward normal derivative on $\Gamma$, $\kappa >0$ is a physical coefficient, 
$\Delta_{\Gamma }$ stands for the Laplace-Beltrami operator
on $\Gamma $ (see, e.g., \cite[Chapter~3]{Gr}), $W_{\Gamma }$ denotes a  potential with some properties similar to those of $W$, and $f_\Gamma$ represents a known datum on $\Sigma$.  

About dynamic boundary conditions, the mathematical research for  
the various problems was already running in 1990's. 
Especially, the Stefan problem with the dynamic boundary condition in the case $\kappa =0$ 
was treated in a series of papers by Aiki~\cite{Ai, Ai2, Ai3}; in particular, 
the existence of a weak solution was investigated. 
Then again, some recent papers dealing with a dynamic boundary condition of type~\eqref{i2} are, among others, \cite{CaCo, CoSp, CGM, GaGu, GaWa, GiMiSc, GoMi, Is, MRSS}.

If one considers the Allen-Cahn equation \eqref{i1} with condition \eqref{i2},
the order parameter $u$ is conserved neither in the bulk nor, as $u_\Gamma$, on the boundary.
The new issue of this paper is the setting of a mass constraint which can involve 
either the solution inside the domain or its counterpart on the boundary (or both of them). 
More precisely, we require that the solution $u$ satisfy 
\begin{equation}
	k_* \le \int_{\Omega }^{}w u(t) dx 
	+ \int_{\Gamma }^{} w_{\Gamma }u_{\Gamma}   (t) d\Gamma 
	\le k^* 
	\quad \mbox{for all }t \in [0,T],
	\label{i3}
\end{equation}
{\pier where} $k_*,k^*$ are given constants fulfilling $k_* \le k^*$, and 
$w$ and $w_{\Gamma }$ are prescribed weight functions on $\Omega $ and $\Gamma$, respectively. 
For example, in the case when $w\equiv 1$, $w_{\Gamma }\equiv 0$ and $k_*=k^*$, 
\eqref{i3} represents the conservation of the volume $\int_{\Omega }^{}{\revis u(x,t)}dx = k_*$, for 
all $t \in [0,T]$, a condition which instead arises naturally from the problem in the framework of a Cahn-Hilliard system (see, e.g., \cite{GiMiSc, Ku}). 

The analysis of the abstract theory for this kind of constraint was 
developed in \cite{FuKe} and motivated from the generalization of concrete problems \cite{Gi, GiSv, SvOm} (see also \cite{Ku}, where the essential structure 
of possible constraints has been discussed for Cahn-Hilliard equation). 
In the abstract approach by \cite{FuKe} the constraint, and in particular the barriers $k_*$ and $k^*$ in \eqref{i3}, are allowed to depend on time. On the other hand, the abstract framework of \cite{FuKe} does not cover a special problem like ours, and especially it does not match with the presence of the {\revis nonlinearities} $W$ and $W_{\Gamma }$. 
In our approach, the solutions of  the system \eqref{i1}--\eqref{i2} are not 
completely free to develop their dynamics, but they should respect the constraint \eqref{i3} on the selected mass values. We discuss and characterize the properties of the unique solution of the initial value problem for the gradient flow system related to \eqref{i1}--\eqref{i3}.

A brief outline of the present paper is as follows. 
In Section~2, we present the main results, consisting in 
the well-posedness of the Allen-Cahn equation with 
dynamic boundary conditions and mass constraints. 
We write the system as an evolution inclusion and characterize the solution
with the help of a Lagrange multiplier. 
In Section~3, we prove the existence result. For the proof, we construct an approximate solution by substituting the maximal monotone graphs with their 
Yosida regularizations. For the approximated problem we 
can apply the result in \cite{FuKe}, by 
checking the validity of the assumptions. Then, after proving some 
uniform estimates, we pass to the limit and conclude the existence proof. 
In Section~4, we prove the continuous dependence: of course, 
this result entails a uniqueness property. 
A final Section~5 contains the proof of a density result, which is useful in our approach.
By the way, here is a detailed index of sections and subsections.}

\begin{itemize}
 \item[1.] Introduction
 \item[2.] Main results
\begin{itemize}
 \item[2.1.] Setting and assumptions
 \item[2.2.] Well-posedness 
 \item[2.3.] Abstract formulations 
\end{itemize}
 \item[3.] Existence
\begin{itemize}
 \item[3.1.] Approximation of the problem
 \item[3.2.] A priori estimates
 \item[3.3.] Passage to the limit
\end{itemize}
 \item[4.] Continuous dependence 
 {\pier \item[5.] Appendix}
\end{itemize}

\section{Main results}
\label{main}
\setcounter{equation}{0}

In this section, we present {\pier our} main result. 
It is the well-posedness of the Allen-Cahn equation with 
dynamic boundary conditions and mass constraints. 
We apply the abstract formulation of the 
evolution inclusion.

\subsection{Setting and assumptions}

Let $0<T<+\infty$ and {\pier let} $\Omega \subset \mathbb{R}^{\takeshi d}$, ${\takeshi d}=2$ or $3$, 
be the bounded domain with smooth boundary $\Gamma:=\partial \Omega $. 
We use the notation{\pier :} 
$H:=L^2(\Omega)$, {\pier $V:=H^1(\Omega)$, 
$H_\Gamma :=L^2(\Gamma)$,} $V_{\Gamma }:=H^1(\Gamma )$, 
with usual norms 
$|\cdot |_H$, $|\cdot |_V$, $|\cdot |_{H_{\Gamma}}$ and $|\cdot |_{V_{\Gamma }}$, 
respectively. 
Then, we obtain
$V \mathop{\hookrightarrow} \mathop{\hookrightarrow} 
H \mathop{\hookrightarrow} \mathop{\hookrightarrow}V^*$, where 
``$\mathop{\hookrightarrow} \mathop{\hookrightarrow} $'' stands for 
the dense and compact embedding, namely 
$(V,H,V^*)$ is {\pier a} standard Hilbert triplet. 
The same considerations {\pier hold for $V_{\Gamma }$ and
$H_{\Gamma }$}. 
Moreover, {\pier we set}
$$
	\mbox{\boldmath $ H $}:=
	H
	\times
	H_{\Gamma },
	\quad 
	\mbox{\boldmath $ V $}:=
	\bigl\{ 
	(u,u_{\Gamma }) \in 
	V \times V_{\Gamma } \, : \  
	u_{|_\Gamma}   =u_{\Gamma } \, \mbox{ a.e.\ on } {\pier \Gamma} 
	\bigr\},
$$
{\pier where $u_{|_\Gamma}   $ denotes the trace of $u$. Observe that 
$\mbox{\boldmath $ H $}$ and $\mbox{\boldmath $ V $}$ are Hilbert spaces} 
with the inner products 
\begin{align*}
	&(\mbox{\boldmath $ u $},\mbox{\boldmath $ z $})_{\mbox{\scriptsize \boldmath $ H$}}
	:=(u,z)_H + (u_\Gamma ,z_\Gamma )_{H_\Gamma } \quad 
	\mbox{for all}~\mbox{\boldmath $ u$}
	:=(u,u_{\Gamma }), \,
	\mbox{\boldmath $ z$}:=(z,z_{\Gamma }) 
	\in \mbox{\boldmath $ H$},\\
	&{\pier (\mbox{\boldmath $ u $},\mbox{\boldmath $ z $})_{\mbox{\scriptsize \boldmath $ V$}}
	:=(u,z)_V + (u_\Gamma ,z_\Gamma )_{V_\Gamma } \quad 
	\mbox{for all}~\mbox{\boldmath $ u$}
	:=(u,u_{\Gamma }), \,
	\mbox{\boldmath $ z$}:=(z,z_{\Gamma }) 
	\in \mbox{\boldmath $ V$}}
\end{align*}
and {\pier related} norms 
\begin{align*}
	&|\mbox{\boldmath $ u$}|_{\mbox{\scriptsize \boldmath $ H$}}
	=\Bigl( |u|_H^2+ |u_{\Gamma }|_{H_{\Gamma }}^2 \Bigr)^{1/2}
	\quad \mbox{for all }\mbox{\boldmath $ u$}:=(u,u_{\Gamma }) \in \mbox{\boldmath $ H$},\\
	&|\mbox{\boldmath $ u$}|_{\mbox{\scriptsize \boldmath $ V$}}
	=\Bigl( |u|_V^2+ \bigl| u_{\Gamma } \bigr|_{V_{\Gamma }}^2 \Bigr)^{1/2}
	\quad \mbox{for all }\mbox{\boldmath $ u$}:=\bigl( u,{\pier u_{\Gamma } }\bigr) \in \mbox{\boldmath $ V$}.
\end{align*}
Then, we obtain 
$\mbox{\boldmath $ V$} 
\mathop{\hookrightarrow} 
\mbox{\boldmath $ H$} 
\mathop{\hookrightarrow} 
\mbox{\boldmath $ V$}^*$, where 
``$\mathop{\hookrightarrow} $'' stands for 
the dense and continuous embedding ({\pier the density is checked in the} Appendix). {\pier
By the way, the above embeddings are also compact, of course.}
As {\pier a remark, let us restate that if $\mbox{\boldmath $ u$}=(u,u_{\Gamma}) \in \mbox{\boldmath $ V$}$ then $u_{\Gamma }$ is exactly the trace of $u$ on $\Gamma$;
while, if $\mbox{\boldmath $ u$}=(u,u_{\Gamma})$ is just in $ \mbox{\boldmath $ H$}$, then} $u \in H$ and $u_{\Gamma } \in H_{\Gamma }$ are independent. 

The {\pier initial-value} problem for the Allen-Cahn equation with 
dynamic boundary conditions is {\pier expressed by} the following {\pier system \eqref{(1)}--\eqref{(3)}} (cf.~\cite{CaCo, CoSp}) 
\begin{gather} 
	\frac{\partial u}{\partial t}-\Delta u + \xi  +\pi(u) = f, \  \hbox{ for some }\   \xi \in \beta(u), 
	\quad \mbox{in } Q, \label{(1)}
\\ 
	u_\Gamma = u_{|_\Gamma}   , \quad 
	{\revis \partial_\nu u
	+ \frac{\partial u_\Gamma}{\partial t}-\Delta_{\Gamma } u_\Gamma
	+\xi_\Gamma +\pi_{\Gamma }(u_\Gamma)  =  f_\Gamma} , \quad  \xi_\Gamma \in 
	\beta_{\Gamma }(u_\Gamma)\quad \mbox{on } \Sigma, \label{(2)}
\\ 
	u(0)=u_0 \quad \mbox{in } \Omega, \quad u_\Gamma(0)=u_{0\Gamma} \quad \mbox{on } \Gamma. 
	\label{(3)}
\end{gather}
where $\beta $, $\beta_{\Gamma }$ are maximal monotone graphs in 
$\mathbb{R} \times \mathbb{R}$. 
{\pier Here, we let $\beta $, $\beta_{\Gamma }$ be the subdifferentials}
$$
	\beta =\partial {\pier\widehat{\beta}}, \quad \beta_{\Gamma }=\partial {\pier\widehat{\beta}}_{\Gamma }
$$
{\pier of some lower semicontinuous and convex functions 
${\pier\widehat{\beta}}${\pier,} ${\pier\widehat{\beta}}_{\Gamma }: 
\mathbb{R} \to [0,+\infty )$ with 
${\pier\widehat{\beta}}(0)={\pier\widehat{\beta}}_{\Gamma }(0)=0$; 
in particular, this implies that $D(\beta )=D(\beta_{\Gamma })=\mathbb{R}$,  
$0 \in \beta (0)$ and $0 \in \beta_{\Gamma }(0)$.} {\pier The given functions
\begin{equation} 
\hbox{$\pi $, $\pi_{\Gamma }: \mathbb{R} \to \mathbb{R}$ are Lipschitz continuous 
with Lipschitz constants $L$, $L_{\Gamma}$,} \label{pilip}
\end{equation}
respectively. Moreover, let 
\begin{equation}
\hbox{$\mbox{\boldmath $ f$}:=(f, f_{\Gamma }) \in L^2(0,T;\mbox{\boldmath $ H$})$ \ \ and 
\ \ {\pier $\mbox{\boldmath $ u$}_0:=(u_0, u_{0\Gamma }) \in \mbox{\boldmath $ V$}$}.}
\label{inidata}
\end{equation}
Now, take an arbitrary  $\mbox{\boldmath $v$} =(v, v_{\Gamma } )$ in, say, $ L^2(0,T;\mbox{\boldmath $V$}) $ and test \eqref{(1)} by $v$: then, with the help of the boundary condition \eqref{(2)} we formally obtain 
{\takeshi 
\begin{align}
& \int_0^T\!\!\!\int_{\Omega}\frac{\partial u}{\partial t} v \, dxdt 
+ \int_0^T\!\!\!\int_{\Omega}\nabla u\cdot\nabla v \, dxdt 
+ \int_0^T\!\!\!\int_{\Omega} \bigl( \xi +\pi(u) \bigr)v dxdt \notag \\
& {}
+ \int_0^T\!\!\!\int_{\Gamma} \frac{\partial u_\Gamma}{\partial t} v_\Gamma \, d\Gamma dt \notag 
+ \int_0^T\!\!\!\int_{\Gamma}\nabla_{\Gamma} u_{\Gamma} \cdot\nabla_{\Gamma}v_{\Gamma} 
\, d\Gamma dt
+ \int_0^T\!\!\!\int_{\Gamma} \bigl( \xi_{\Gamma} +\pi_{\Gamma}(u_{\Gamma}) \bigr)v_\Gamma 
\, d\Gamma dt \notag \\
& = \int_0^T\!\!\!\int_{\Omega}f v \, dxdt+ \int_0^T\!\!\!\int_{\Gamma}f_{\Gamma} v_\Gamma 
\, d\Gamma dt, 
\label{p1}
\end{align}
}a  variational equality holding for all 
$\mbox{\boldmath $v$} =(v, v_{\Gamma }) \in L^2(0,T;\mbox{\boldmath $V$}) $
and yielding a weak formulation of \eqref{(1)}--\eqref{(2)}. Here, $\nabla_\Gamma$ denotes the surface gradient on $\Gamma$ (see,\ e.g.,\ \cite[Chapter 3]{Gr}). Hence, we can argue that a suitable space for the solution $\mbox{\boldmath $u$} =(u, u_{\Gamma } )$ of \eqref{(3)}--\eqref{p1} (in which $\xi$ and $\xi_\Gamma$ represent selections of  $\beta(u) $ and $\beta_{\Gamma }(u_\Gamma) $ as in is \eqref{(1)}--\eqref{(2)}) is 
$H^1(0,T;\mbox{\boldmath $H$}) \cap L^2(0,T; \mbox{\boldmath $V$})$. Under suitable conditions on 
$\beta, \, \beta_\Gamma, \, u_0, u_{0\Gamma }$ such a (unique) solution actually exists (see, e.g., \cite{CaCo}) and possesses further regularity properties: $\mbox{\boldmath $\xi$} =(\xi, \xi_{\Gamma } )\in L^2(0,T; \mbox{\boldmath $H$}) $ and, in particular, 
$\mbox{\boldmath $ u$}:=(u, u_{\Gamma }) \in L^2(0,T;H^2(\Omega) ) \times L^2(0,T;H^2(\Gamma) ) $.

On the other hand, in this paper we are interested to the variational inequality 
obtained by \eqref{p1} when replacing the ``$=$'' sign by ``$\leq$'' and taking, in place of the
test element $\mbox{\boldmath $v$}$, the difference 
$\mbox{\boldmath $u$} - \mbox{\boldmath $z$}$, 
where both the solution $\mbox{\boldmath $u$} =(u, u_{\Gamma } )$ and the arbitrary 
$\mbox{\boldmath $z$} =(z, z_{\Gamma }) \in L^2(0,T;\mbox{\boldmath $V$})$
have to satisfy the constraint (written in terms of $z$ and $z_\Gamma$)
\begin{equation} 
	k_* \le \int_{\Omega }^{}w \, z(t) dx 
	+ \int_{\Gamma }^{} w_{\Gamma }\, z_\Gamma (t) d\Gamma 
	\le k^* ,
	\quad t \in [0,T]. \label{(4)}
\end{equation}
Here, $k_*$ and $k^*$ are real constants with $k_* \le k^*$, and 
$\mbox{\boldmath $ w$}:=(w,w_{\Gamma }) $ is fixed in $\mbox{\boldmath $ H$}$. 
We require that the weight functions $w$ and $w_{\Gamma }$ satisfy
$$
	w \ge 0 
	\quad 
	\mbox{a.e.\ in } \Omega, \quad 
	w_{\Gamma } \ge 0 
	\quad 
	\mbox{a.e.\ on } \Gamma
$$
and
\begin{equation}
	{\takeshi \sigma_0} :=\int_{\Omega}^{} w dx + \int_{\Gamma }^{} w_{\Gamma } d\Gamma >0.
\label{p2}
\end{equation}
The constraint \eqref{(4)} entails a limitation on a specific averaged value of the solution $\mbox{\boldmath $u$} =(u, u_{\Gamma } )$ in the bulk and/or on the boundary. 
Inequality \eqref{p2} can be seen as a non-degeneracy condition on the weight element 
$\mbox{\boldmath $ w$}:=(w,w_{\Gamma })$.}

{\pier Hence, let us term (P) the initial-value problem related to the variational inequality and
to the constraint in  \eqref{(4)}. Now, we define precisely the notion of solution to the problem (P) by means of a Lagrange multiplier.} 

\paragraph{Definition 2.1.} 
{\pier {\it The triplet $(\mbox{\boldmath $ u$}, \mbox{\boldmath $ \xi$}, \lambda )$  is called the solution of {\rm (P)} if
\begin{align*}
	\mbox{{\boldmath $ u$}}=(u,u_{\Gamma }) \quad \hbox{with} \quad
	u \in H^1(0,T;H) \cap C ([0,T];V) \cap L^2 {\revis \bigl(} 0,T;H^2(\Omega ) {\revis \bigr)},\\
	u_{\Gamma } \in H^1(0,T;H_{\Gamma}) \cap C([0,T];V_{\Gamma }) \cap L^2 {\revis \bigl(} 0,T;H^2(\Gamma ){\revis \bigr)},\\
	\mbox{{\boldmath $\xi$}}=(\xi ,\xi_{\Gamma }) \in  L^2(0,T;\mbox{\boldmath $H$}) , \quad
	\lambda \in L^2(0,T) 
\end{align*}
and $u,\, u_\Gamma, \, \xi,\, \xi_\Gamma, \,\lambda$ satisfy
\begin{gather} 
	\frac{\partial u}{\partial t}-\Delta u + \xi + \pi(u) + \lambda w = f
	\quad \mbox{a.e.\ in } Q, \label{(8)}
\\
	\xi \in \beta (u)
	\quad \mbox{a.e.\ in } Q,\label{(9)}
\\ 
	u_{|_\Gamma}    =u_\Gamma, \quad 
	\partial_\nu u 
	+ \frac{\partial u_{\Gamma }}{\partial t}-\Delta_{\Gamma } u_{\Gamma }
	+ \xi_{\Gamma } +\pi_{\Gamma }(u_{\Gamma }) + \lambda w _{\Gamma } = f_\Gamma 
	\quad \mbox{a.e.\ on } \Sigma, \label{(10)}
\\ 
	\xi_\Gamma \in \beta_{\Gamma } (u_{\Gamma })
	\quad \mbox{a.e.\ on } \Sigma, \label{(11)}
\\
	u(0)=u_0 \quad \mbox{a.e.\ in } \Omega, \quad 
	u_{\Gamma }(0)=u_{0\Gamma} \quad \mbox{a.e.\ on } \Gamma, \label{(12)}
\\ 
	k_* \le \int_{\Omega}^{} w u(t) dx + \int_{\Gamma }^{} w_{\Gamma } u_{\Gamma }(t)d\Gamma 
	\le k^*
	\quad \mbox{for all\ } t \in [0,T], \label{(13)}
\\
\lambda (t) \left(\int_{\Omega}^{} w 
{\fukao \bigl( }u(t)- z {\fukao \bigr)}  dx 
+ \int_{\Gamma }^{} w_{\Gamma } {\fukao \bigl( }u_{\Gamma }(t) - z_\Gamma 
{\fukao \bigr)} d\Gamma \right) \geq 0 \nonumber \quad  \hbox{for a.a. } t \in (0,T)\qquad\\ 
\mbox{and for all\ }  \mbox{\boldmath $z$} =(z, z_{\Gamma }) 
\in \mbox{\boldmath $V$}
\  \mbox{such that } \ k_* \le \int_{\Omega}^{} w z dx + 
\int_{\Gamma }^{} w_{\Gamma } z_{\Gamma }d\Gamma \le k^*.
 \label{13bis}
\end{gather} 
}}

\subsection{Well-posedness} 

{\pier The first result states the continuous dependence on the data. 
The uniqueness of the component $\mbox{{\boldmath $ u$}}$ of the solution 
(see the later Remark~3.3) is also guaranteed by this theorem. 

\paragraph{Theorem 2.1.}
{\it
For $i=1,2$ let $(\mbox{\boldmath $ u$}^{(i)}, \mbox{\boldmath $ \xi$}^{(i)}, \lambda^{(i)} )$,
with $\mbox{\boldmath $ u$}^{(i)}= (u^{(i)}, u^{(i)}_{\Gamma })$ 
and $\mbox{\boldmath $ \xi$}^{(i)}= (\xi^{(i)}, \xi^{(i)}_{\Gamma })$,
be a solution to {\rm (P)} corresponding to the data $\mbox{\boldmath $ f$}^{(i)}= (f^{(i)}, f^{(i)}_{\Gamma })$ and $ \mbox{\boldmath $ u$}_0^{(i)} = (u_0^{(i)}, u^{(i)}_{0\Gamma })$. 
Then, there exists a positive constant $C>0$, depending only on  $L$, $L_{\Gamma}$ and $T$, such that 
\begin{eqnarray} 
	\lefteqn{ 
	\bigl| u^{(1)}(t)-  u^{(2)}(t) \bigr|_H^2 
	+
	\bigl| u^{(1)}_{\Gamma }(t)-  u^{(2)}_{\Gamma }(t) \bigr|_{H_\Gamma }^2 
	} \nonumber \\
	& ~ & {} + 2 \int_{0}^{t} 
	\bigl| \nabla u^{(1)}(s) -\nabla u^{(2)}(s ) \bigr|_{H^{\takeshi d}}^2 ds 
	+ 2 \int_{0}^{t} 
	\bigl| \nabla _\Gamma u^{(1)}_{\Gamma }(s ) -\nabla _\Gamma u^{(2)}_{\Gamma }(s ) 
	\bigr|_{H_{\Gamma }^{{\takeshi d}}}^2 ds \nonumber
	\\
	&\le &C 
	\left\{ 
	\bigl| u^{(1)}_0-  u^{(2)}_0 \bigr|_H^2 
	+
	\bigl| u^{(1)}_{0\Gamma }-  u^{(2)}_{0\Gamma } \bigr|_{H_\Gamma }^2 
	+ \int_{0}^{T} 
	\bigl| f^{(1)}(s ) -f^{(2)}(s ) \bigr|_H^2 ds 
	\right. \nonumber \\
	& ~ & \left. \qquad {}
	+ \int_{0}^{T} 
	\bigl| f^{(1)}_{\Gamma }(s ) -f^{(2)}_{\Gamma }(s ) \bigr|_{H_{\Gamma }}^2 ds 
	\right\}  \quad \hbox{for all } \, t \in [0,T]. \label{p0}
\end{eqnarray} 
}

The second result deals with the existence of the solution. 
To this aim, we further assume that there exist positive constants 
{\pier $c_0,\, {\takeshi \varrho}$} such that 
\begin{gather} 
	|s| \le c_0 \bigl( 1+{\pier\widehat{\beta}}(r) \bigr) \quad 
	\mbox{for all } r \in \mathbb{R}~
	\mbox{and } s \in \beta (r),\label{(5)}
\\
	|s| \le c_0 \bigl( 1+{\pier\widehat{\beta}}_{\Gamma }(r) \bigr) \quad 
	\mbox{for all } r \in \mathbb{R}
	\mbox{ and } s \in \beta_{\Gamma }(r),\label{(6)}
\\
	{\fukao \bigl|}\beta^\circ (r) {\fukao \bigr|} 
	\le {\takeshi \varrho} {\fukao \bigl|}\beta_{\Gamma }^\circ (r){\fukao \bigr|}+c_0
	\quad 
	\mbox{for all } r \in \mathbb{R},\label{(7)}
\end{gather} 
{\pier where} the minimal section $\beta^\circ $ of $\beta $ is {\pier specified} by
$$
	\beta^\circ (r):=\bigl\{ r^* \in \beta (r)\, : \  |r^*|=\min _{s \in \beta (r)} |s| \bigr\}, \quad r\in \mathbb{R}
$$
and the same definition holds for $\beta_{\Gamma }^\circ $ (and for any maximal monotone graph!). We also require compatibility conditions for the initial data, that are
\begin{equation}
	k_* \le \int_{\Omega}^{} w u_0 dx + \int_{\Gamma }^{} w_{\Gamma } u_{0\Gamma }d\Gamma \le k^*  \label{p3}
\end{equation}
and 
\begin{equation}
{\pier {\pier\widehat{\beta}} (u_0) \in L^1(\Omega ), \quad  {\pier\widehat{\beta}}_{\Gamma }(u_{0\Gamma}) \in L^1(\Gamma).} \label{p4}
\end{equation}

\paragraph{Theorem 2.2.} 
{\it 
Under the above assumptions, there exists one solution of {\rm (P)}. 
}}

\subsection{Abstract formulation} 

{\pier In this subsection, we comment on the formulation of the problem and on our results. 
The first remark} is related to the mathematical treatment by the evolution inclusion 
governed by subdifferential operators. 

Our mass constraint \eqref{(13)} {\pier can be rewritten as} 
$$
	k_* \le \bigl( \mbox{\boldmath $ w$},\mbox{\boldmath $ u$}(t)
	\bigr)_{\mbox{\scriptsize \boldmath $ H$}} \le k^* 
	\quad 
	\mbox{for all }
	t \in [0,T].
$$
Then, we define the convex constraint set $\mbox{\boldmath $ K$}$ 
{\pier that} plays {\pier an} important role in this paper:
$$ 
	\mbox{\boldmath $ K$}:=
	\bigl\{ \mbox{\boldmath $ z$} \in \mbox{\boldmath $ V$} \, : \ 
	k_* \le ( \mbox{\boldmath $ w$},\mbox{\boldmath $ z$}
	)_{\mbox{\scriptsize \boldmath $ H$}} \le k^* 
	\bigr\},
$$
with the indicator function $I_{\mbox{\boldmath \scriptsize $ K$}}:\mbox{\boldmath $ H$} \to [0,+\infty ]$ {\pier fulfilling $I_{\mbox{\boldmath \scriptsize $K$}} (\mbox{\boldmath $ z$}) = 0 $ if $ \mbox{\boldmath $ z$} \in \mbox{\boldmath $ K$}$, $I_{\mbox{\boldmath \scriptsize $K$}} (\mbox{\boldmath $ z$}) = +\infty $ if $	\mbox{\boldmath $ z$} \in \mbox{\boldmath $ H$}\setminus \mbox{\boldmath $ K$}. $
Moreover, we {\pier introduce the} proper, lower semicontinuous and convex functional 
$\varphi: \mbox{\boldmath $ H$} \to [0,+\infty ]$ by 
$$
	\varphi (\mbox{\boldmath $ z$}) 
	:= \left\{ 
	\begin{array}{l}
	\displaystyle 
	\frac{1}{2} \int_{\Omega }^{} \bigl |\nabla z \bigr |^2 dx + \int_{\Omega }^{} 
	{\pier\widehat{\beta}}(z)dx + \frac{1}{2} \int_{\Gamma }^{} |\nabla _{\Gamma } z_{\Gamma }|^2 d\Gamma  
	+ \int_{\Gamma }^{} 
	{\pier\widehat{\beta}}_{\Gamma } (z_{\Gamma })d\Gamma \vspace{2mm}\\
	\hfill
	\mbox{if } 
	\mbox{\boldmath $ z$} \in \mbox{\boldmath $ V$}, 
	{\pier\widehat{\beta}}(z) \in L^1(\Omega )~
	\mbox{and } {\pier\widehat{\beta}}_{\Gamma }(z_{\Gamma }) \in L^1(\Gamma), \vspace{2mm}\\
	+\infty \quad \mbox{otherwise} . 
	\end{array} 
	\right. 
$$
Therefore, it is possible to check that our problem enters the following abstract 
form of an evolution inclusion with a Lipschitz perturbation:}
\begin{equation} 
	\mbox{\boldmath $ u$}'(t)+\partial (\varphi +I_{\mbox{\boldmath \scriptsize $ K$}})
	\bigl( \mbox{\boldmath $ u$}(t) \bigr)
	+ \mbox{\boldmath $ \pi $} \bigl (\mbox{\boldmath $ u$}(t) \bigr )
	\ni \mbox{\boldmath $ f$}(t) 
	\quad \mbox{in } \mbox{\boldmath $ H$}, 
	\, \mbox{ for a.a.\ } t \in (0,T), \label{(14)}
\end{equation} 
where
$\mbox{\boldmath $ \pi $}(\mbox{\boldmath $ z$})
:=( \pi (z),\pi_{\Gamma }(z_\Gamma ))$ for all $\mbox{\boldmath $ z$} \in \mbox{\boldmath $ H$}$. 
This kind of evolution inclusion is {\colli well known} as a gradient flow equation 
including a Lipschitz perturbation, and it has been treated, in particular, in \cite{Br}. 

{\pier Thus, from this point of view the existence and uniqueness of the solution to the Cauchy problem 
for \eqref{(14)} is perfectly known. On the other hand, what is important here is to characterize the 
suitable selection from  $\partial (\varphi +I_{\mbox{\boldmath \scriptsize $ K$}})
	\bigl( \mbox{\boldmath $ u$}(t) \bigr)$ for a.a.\ $t\in (0,T)$, which is our main concern. 
Now, one can check that (see, e.g., \cite[p.\ 59]{Ba} or \cite{CaCo}) the subdifferential operator $\partial \varphi $ can be expressed in a formal way as 
\begin{align*}
&\mbox{\boldmath $ z$}^*:=(z^*,z^*_{\Gamma }) 
\in \partial \varphi (\mbox{\boldmath $ z$}) \hbox{ is in 
\mbox{\boldmath $ H$} if and~only~if}
\\
&(z^*,z^*_{\Gamma }) = \bigl( -\Delta z+\beta (z), \partial_\nu z-\Delta_{\Gamma } z_\Gamma  + \beta_{\Gamma }(z_\Gamma ) \bigr).
\end{align*}
{\takeshi However}, when one adds the indicator function  $I_{\mbox{\boldmath \scriptsize $ K$}}$ to $\varphi$,
then the subdifferential $\partial (\varphi +I_{\mbox{\boldmath \scriptsize $ K$}})$ must take 
into account the constraint given by $\mbox{\boldmath $ K$}$. Then, the point of emphasis of Theorem 2.2 is the (further) characterization with the 
help of the Lagrange multiplier $\lambda $ in \eqref{(8)} and \eqref{(10)}. Indeed, our analysis shows in particular that 
$ \mbox{\boldmath $ z $}^*:=(z^*,z^*_{\Gamma }) \in \partial (\varphi + I_{\mbox{\boldmath 
\scriptsize $ K$}}) (\mbox{\boldmath $ z$})$ lies in $ \mbox{\boldmath $ H$}$ if and only if there is a scalar $\lambda_{\mbox{\boldmath \scriptsize $ z$}}$ such that
$$
	(z^*,z^*_{\Gamma }) = \bigl( -\Delta z+\beta (z) + \lambda_{\mbox{\boldmath \scriptsize $ z$}} w, \partial_\nu z-\Delta_{\Gamma } z_\Gamma  + \beta_{\Gamma }(z_\Gamma ) + \lambda_{\mbox{\boldmath \scriptsize $ z$}} w_\Gamma \bigr).
$$
We point out that such a characterization problem was already treated in \cite{FuKe}, 
in an abstract framework with appropriate assumptions for the abstract functions and tools. 
However, in our concrete problem for the Allen-Cahn equation with dynamic boundary conditions 
we cannot ensure the validity of \cite[Assumption~(A2)]{FuKe} for $\varphi $. 
Therefore, in the next section we consider an approximating problem to which we can (more or less) apply the abstract result of \cite{FuKe}. In this sense, our results turn out to be an extension of \cite{FuKe} for our concrete problem.}

\section{Existence}
\setcounter{equation}{0}

{\pier This section is devoted to the proof of {\pier Theorem~2.2.} 
We make use of {\pier the Yosida approximations for maximal monotone operators and of well-known 
results of this theory (see, e.g., \cite{Ba,Br,Ke}).} For each $\varepsilon \in (0,1]$, we define 
$\beta_\varepsilon, \beta_{\Gamma,\varepsilon}:\mathbb{R} \to \mathbb{R}$, along with
the associated resolvent operators $J_\varepsilon, J_{\Gamma,\varepsilon}:\mathbb{R} \to \mathbb{R}$, by 
$$
	\beta_\varepsilon (r)
	:= \frac{1}{\varepsilon } {\revis \bigl(} r-J_\varepsilon (r) {\revis \bigr)}
	:=\frac{1}{\varepsilon }\bigl( r-(I+\varepsilon \beta )^{-1} (r)\bigr),
$$
$$
	\beta_{\Gamma, \varepsilon} (r)
	:= \frac{1}{\varepsilon {\takeshi \varrho} } {\revis \bigl(} r-J_{\Gamma,\varepsilon  }(r) {\revis \bigr)}
	:=\frac{1}{\varepsilon {\takeshi \varrho}  }\bigl( r- (I+\varepsilon {\takeshi \varrho}  \beta_\Gamma )^{-1} (r) \bigr)
	\quad \mbox{for all } r \in \mathbb{R},
$$
where ${\takeshi \varrho}  >0$ is same constant as in \eqref{(7)}. Note that the two definitions are not symmetric since in the second it is ${\takeshi \varrho}  \varepsilon $ and not directly $ \varepsilon $ to be used as approximation parameter. Anyway, we easily have 
$\beta_\varepsilon(0)=\beta_{\Gamma, \varepsilon}(0)=0$. 
Moreover, the related Moreau-Yosida regularizations ${\pier\widehat{\beta}}_\varepsilon, \, {\pier\widehat{\beta}}_{\Gamma,\varepsilon}$
of ${\pier\widehat{\beta}}, \, {\pier\widehat{\beta}}_{\Gamma } : \mathbb{R} \to \mathbb{R}$
fulfill
\begin{align*}
	{\pier\widehat{\beta}}_{\varepsilon }(r)
	:=\inf_{s \in \mathbb{R}}\left\{ \frac{1}{2\varepsilon } |r-s|^2+{\pier\widehat{\beta}}(s) \right\} 
	= 
	\frac{1}{2\varepsilon } {\revis \bigl|}r-J_\varepsilon (r){\revis \bigr|}^2
	+{\pier\widehat{\beta}}{\revis \bigl(}J_\varepsilon( r){\revis \bigr)}
	= \int_{0}^{r} \beta_\varepsilon (s)ds,
\\
	{\pier\widehat{\beta}}_{\Gamma, \varepsilon }(r)
	:=\inf_{s \in \mathbb{R}}\left\{ \frac{1}{2\varepsilon {\takeshi \varrho}  } |r-s|^2+{\pier\widehat{\beta}}_\Gamma (s) \right\} 
	= \int_{0}^{r} \beta_{\Gamma,\varepsilon} (s)ds
	\quad \mbox{ for all } r\in \mathbb{R}.
\end{align*}
It is {\fukao well known} that $\beta_\varepsilon$ is Lipschitz continuous with Lipschitz constant 
$1/\varepsilon $ and $\beta_{\Gamma, \varepsilon}$ is also Lipschitz continuous with constant  $1/(\varepsilon {\takeshi \varrho} )$. In addition, we have the standard properties 
$$
	\bigl |\beta_\varepsilon (r) \bigr | \le \bigl |\beta^\circ (r) \bigr |, \quad 
	\bigl |\beta_{\Gamma ,\varepsilon }(r) \bigr | \le \bigl |\beta_{\Gamma }^\circ (r) \bigr | 
	\quad \mbox{for all } r \in \mathbb{R},
$$ 
$$
	0 \le \widehat{\beta}_\varepsilon (r) \le \widehat{\beta}(r), \quad 
	0 \le \widehat{\beta}_{\Gamma, \varepsilon} (r) \le \widehat{\beta}_{\Gamma }(r)
	\quad \mbox{for all } r \in \mathbb{R}.
$$

We emphasize that \eqref{(5)}--\eqref{(6)} entail
\begin{equation} 
	\bigl |\beta_\varepsilon (r) \bigr | \le c_0 
	\bigl( 1+\widehat{\beta}_\varepsilon (r) \bigr) \quad 
	\mbox{for all } r \in \mathbb{R}, \label{(15)}
\end{equation}
\begin{equation} 
	\bigl |\beta_{\Gamma,\varepsilon }(r) \bigr | \le c_0 \bigl( 1+\widehat{\beta}_{\Gamma, \varepsilon  }(r) \bigr) \quad 
	\mbox{for all } r \in \mathbb{R},\label{(16)}
\end{equation}
with the same constant $c_0$ as in \eqref{(5)}--\eqref{(6)}. Indeed, arguing for instance for $\beta_\varepsilon$, it suffices to notice that for all $r \in \mathbb{R}$ there exists $s_\varepsilon \in \beta( J_\varepsilon (r)) $ such that 
$$
	\bigl |\beta_\varepsilon (r) \bigr | = |s_\varepsilon| \leq c_0 
	{\revis \Bigl(} 1+ \widehat{\beta} {\revis \bigl(} J_\varepsilon (r) {\revis \bigr)} {\revis \Bigr)} = c_0 
	{\revis \Bigl(} 1+ \widehat{\beta}_\varepsilon (r){\revis \bigr)} {\revis \Bigr)},
$$
thanks to \eqref{(5)}. Moreover, owing to the assumption~\eqref{(7)} and \cite[Lemma 4.4]{CaCo}, 
the inequality
\begin{equation}
	\bigl |\beta_\varepsilon (r)\bigr | \le {\takeshi \varrho}  \bigl |\beta_{\Gamma,\varepsilon} (r)\bigr |+c_0
	\quad 
	\mbox{for all } r \in \mathbb{R},\label{(17)}
\end{equation} 
holds for $\beta_\varepsilon$ and $\beta_{\Gamma,\varepsilon}$ as well.}

\subsection{Approximation of the problem}

{\pier Let us consider an approximation of (P) which is stated as the following initial-value problem for a gradient flow equation: for each $\varepsilon \in (0,1]$ let $\mbox{\boldmath $ u$}_\varepsilon$ solve the abstract Cauchy problem}
\begin{gather} 
	\mbox{\boldmath $ u$}_\varepsilon '(t)+\partial (\varphi _\varepsilon +I_{\mbox{\boldmath \scriptsize $ K$}})
	\bigl( \mbox{\boldmath $ u$}_\varepsilon (t) \bigr) 
	+\mbox{\boldmath $ \pi $}\bigl( \mbox{\boldmath $ u$}_\varepsilon (t) \bigr)
	\ni \mbox{\boldmath $ f$}(t)
	\quad \mbox{in } \mbox{\boldmath $ H$}{\pier ,} 
	\, \mbox{ for a.a.\ } t \in (0,T),\label{(18)}
\\ 
	\mbox{\boldmath $ u$}_\varepsilon (0)=\mbox{\boldmath $ u$}_0 
	\quad \mbox{in } \mbox{\boldmath $ H$},\label{(19)}
\end{gather} 
with $\mbox{\boldmath $ u$}_0=( u_0,u_{0\Gamma} ) \in \mbox{\boldmath $ K$}$ {\pier satisfying \eqref{p4}
and  $\varphi_\varepsilon:\mbox{\boldmath $ H$} \to [0, +\infty ]$ being} defined by
$$ 
	\varphi_\varepsilon (\mbox{\boldmath $ z$}) 
	:=
	\left\{ 
	\begin{array}{l} 
	\displaystyle 
	\frac{1}{2} \int_{\Omega }^{} |\nabla z |^2 dx 
	+ \int_{\Omega }^{} 
	\widehat{\beta}_\varepsilon ( z )dx 
	+\frac{\varepsilon }{2} \int_{\Omega }^{} |z|^2dx \vspace{2mm}\\
	\displaystyle \quad
	{} + \frac{1}{2} \int_{\Gamma }^{} |\nabla _{\Gamma } z_{\Gamma }|^2 d\Gamma 
	+ \int_{\Gamma }^{} 
	\widehat{\beta}_{\Gamma,\varepsilon } ( z_{\Gamma }) d\Gamma 
	+ \frac{\varepsilon }{2} \int_{{\pier \Gamma} }^{} |z_{\Gamma}|^2{\pier d\Gamma} \quad 
	\mbox{if } \mbox{\boldmath $ z$} \in \mbox{\boldmath $ V$}, \vspace{2mm}\\
	+\infty \quad \mbox{if } \mbox{\boldmath $ z$} \in \mbox{\boldmath $ H$} \setminus \mbox{\boldmath $ V$};
	\end{array} 
	\right. 
$$
moreover, it is understood that
$$\mbox{\boldmath $ \pi $}(\mbox{\boldmath $ z$})
:={\fukao \bigl( } \pi (z),\pi_{\Gamma }(z_\Gamma ) {\fukao \bigr)} 
\quad \hbox{for all } \mbox{\boldmath $ z$}= ( z, z_\Gamma) \in \mbox{\boldmath $ H$}.$$ 

For a proper convex lower semicontinuous function $\psi :\mbox{\boldmath $ V$} 
\to (-\infty , +\infty ]$, we denote
by $\partial_* \psi $ its subdifferential operator acting from  $\mbox{\boldmath 
$ V$}$ to $\mbox{\boldmath $ V$}^*$. In the next statement 
we point out the following characterization 
of $\partial_* \varphi_\varepsilon \, $.

\paragraph{Lemma 3.1.} 
{\it The function $\varphi_\varepsilon:\mbox{\boldmath $ H$} \to [0, +\infty ]$ is convex and lower semicontinuous, with domain $D(\varphi_\varepsilon )= \mbox{\boldmath $ V$}$. Moreover, 
$\varphi_\varepsilon$ is lower semicontinuous in $\mbox{\boldmath $ V$}$ as well and 
the subdifferential 
$\partial_* \varphi _\varepsilon $ is single-valued and specified by the following form:
\begin{eqnarray} 
	\bigl \langle 
	\partial_* \varphi _\varepsilon(\mbox{\boldmath $ z$}), 
	\bar{\mbox{\boldmath $ z$}}
	\bigr \rangle_{\mbox{\boldmath \scriptsize $ V$}^*,\mbox{\boldmath \scriptsize $ V$}} 
	& = & \bigl (\nabla z,\nabla \bar{z} \bigr )_{H^{{\takeshi d}}}+ 
	\bigl (\beta_\varepsilon (z),\bar{z} \bigr )_{H}+
	\varepsilon \bigl (z, \bar{z} \bigr )_{H}
	\nonumber \\
	& ~ & 
	{}+
	\bigl (\nabla _{\Gamma } z_{\Gamma },\nabla _{\Gamma } \bar{z}_{\Gamma } \bigr )_{H_{\Gamma}^{{\takeshi d}} }
	+
	\bigl (\beta_{\Gamma ,\varepsilon }(z_{\Gamma }),\bar{z}_{\Gamma } \bigr )_{H_{\Gamma }}
	+\varepsilon \bigl (z_\Gamma , \bar{z}_\Gamma \bigr )_{H_\Gamma }\nonumber \\
	& ~ & 
	{}\qquad \mbox{ for all } 
	\mbox{\boldmath $ z$}=(z,z_{\Gamma }),\,  
	\bar{\mbox{\boldmath $ z$}}=(\bar{z},\bar{z}_{\Gamma }) 
	\in \mbox{\boldmath $ V$}. \label{(20)}
\end{eqnarray} 
Finally, there exists a positive constant 
$C_\varepsilon >0$ depending on $\varepsilon >0$ such that
\begin{equation} 
	\bigl |
	\partial_* \varphi _\varepsilon (\mbox{\boldmath $ z$}) 
	\bigr |_{\mbox{\boldmath \scriptsize $ V$}^*} 
	\le C_\varepsilon \bigl(1+\varphi _\varepsilon (\mbox{\boldmath $ z$}) \bigr) \quad 
	\mbox{for all } \mbox{\boldmath $ z$} \in \mbox{\boldmath $ V$}.\label{(21)}
\end{equation} 
}%

\paragraph{Proof.} 
The function $\varphi_\varepsilon$ is convex and 
assumes finite value on all elements of $\mbox{\boldmath $ V$}$; in addition, 
it is straightforward to check that $\varphi_\varepsilon$ is strongly, 
whence also weakly, lower semicontinuous in $\mbox{\boldmath $ V$}$. 
Now, let $\mbox{\boldmath $ z$}_n \to \mbox{\boldmath $ z$} $ 
strongly in $\mbox{\boldmath $ H$}$ as $n\to \infty$, 
and assume that $ \varphi_\varepsilon (\mbox{\boldmath $ z$}_n ) \leq \alpha $ 
for some $\alpha\geq 0 $ and all $n \in \mathbb{N}$. 
Then, as $\widehat{\beta}_{\varepsilon }$ and $\widehat{\beta}_{\Gamma ,\varepsilon } $ 
are non-negative, 
we easily conclude that $\{ \mbox{\boldmath $ z$}_n {\revis \}_{n \in \mathbb{N}}}$ is 
bounded in $\mbox{\boldmath $ V $}$ and consequently ${\colli \mbox{\boldmath $ z$}_n}$  weakly converges 
to $\mbox{\boldmath $ z$} $ in $\mbox{\boldmath $ V $}$ as $n\to \infty$. 
Then, it turns out that 
$ \varphi_\varepsilon (\mbox{\boldmath $ z$}) \leq \alpha $ 
and the weak lower semicontinuity of $ \varphi_\varepsilon$ in $\mbox{\boldmath $ H$}$ follows. 
Next, let  $\mbox{\boldmath $ z$}^* \in \partial_* \varphi _\varepsilon 
(\mbox{\boldmath $ z$})$ in $\mbox{\boldmath $ V$}^*$. 
Then, from the definition of the subdifferential, we have 
\begin{eqnarray} 
	\lefteqn{ 
	\langle \mbox{\boldmath $ z$}^*, 
	\delta \bar{\mbox{\boldmath $ z$}} \rangle_{\mbox{\boldmath \scriptsize $ V$}^*,
	\mbox{\boldmath \scriptsize $ V$}}
	} \nonumber \\  
	& \le & 
	\delta \int_{\Omega }^{} \nabla z \cdot \nabla \bar{z} dx 
	+ \frac{\delta^2}{2} \int_{\Omega }^{} |\nabla \bar{z}|^2 dx 
	+ \int_{\Omega }^{} 
	\left\{ \widehat{\beta }_\varepsilon (
	z+\delta \bar{z})-\widehat{\beta}_\varepsilon ( z ) \right\} dx 
	\nonumber \\
	& ~ & {} 
	+\delta \varepsilon \int_{\Omega }^{} z \bar{z} dx 
	+ \frac{\delta ^2 \varepsilon }{2} \int_{\Omega }^{} 
	|\bar{z}|^2dx 
	+\delta \int_{\Gamma }^{} \nabla _{\Gamma} z_{{\fukao \Gamma }} \cdot \nabla_{\Gamma } \bar{z}_{{\fukao \Gamma }} d\Gamma 
	+ \frac{\delta^2}{2} \int_{\Gamma }^{} 
	|\nabla_{\Gamma } 
	\bar{z}_{\Gamma }|^2  d\Gamma  
	\nonumber \\
	& ~ & {}
	 + \int_{\Gamma }^{} 
	 \left\{ \widehat{\beta }_{\Gamma,\varepsilon}(
	z_{\Gamma }
	+\delta \bar{z}_{\Gamma })-\widehat{\beta}_{ \Gamma, \varepsilon}(z_{{\fukao \Gamma }})
	\right\} d\Gamma
	+\delta \varepsilon \int_{\Gamma }^{} z_\Gamma \bar{z}_\Gamma d\Gamma
	+ \frac{\delta ^2\varepsilon }{2} \int_{\Gamma }^{} |\bar{z}_\Gamma |^2 d\Gamma,
	\label{(22)} 
\end{eqnarray}
for all $\bar{\mbox{\boldmath $ z$}} \in \mbox{\boldmath $ V$}$ 
and $\delta >0$. 
Here, from the Lipschitz continuity of $\beta_\varepsilon $ and ${\beta}_{\Gamma,\varepsilon }$
we infer that
\begin{align*}
&	\left| \frac{\, 
	\widehat{\beta}_\varepsilon (z + \delta \bar{z})-\widehat{\beta}_\varepsilon (z)}{\delta } \right| 
	\le {\fukao \bigl|} {\beta}_\varepsilon (\zeta)-\beta_\varepsilon (0) {\fukao \bigr|}
	|\bar{z}|
	\leq
	\frac{1}{\varepsilon }{\fukao \bigl(}
	|z| + \delta|\bar{z}|{\fukao \bigr)} 
	|\bar{z}| \quad \mbox{a.e.\ in } \Omega,  \\
& \left| \frac{	\, 
 \widehat{\beta }_{\Gamma,\varepsilon}( z_{\Gamma }
	+\delta \bar{z}_{\Gamma })-\widehat{\beta}_{\Gamma,\varepsilon }(z_{{\fukao \Gamma }})}\delta
	\right|	\le {\fukao \bigl|}
	{\beta}_{\Gamma,\varepsilon} (\zeta_{\Gamma })-
	\beta_{\Gamma,\varepsilon} (0) {\fukao \bigr|}
	|\bar{z}_{\Gamma }|	\leq
	\frac{1}{\varepsilon {\takeshi \varrho}  } 
	{\fukao \bigl(}|z_{\Gamma }| + \delta|\bar{z}_{\Gamma }|{\fukao \bigr)} |\bar{z}_{\Gamma }|
	 \quad \mbox{a.e.\ on } \Gamma, 	
\end{align*}
for some intermediate functions $\zeta :\Omega \to \mathbb{R} $, between $z$ and $\bar{z}$, and $\zeta_{\Gamma }:\Gamma \to \mathbb{R}$, between $z_{\Gamma }$ and $\bar{z}_{\Gamma }$.
 Therefore, dividing \eqref{(22)} by 
$\delta $ and letting $\delta \to 0$, we obtain 
\begin{eqnarray*} 
	\bigl \langle 
	\mbox{\boldmath $ z$}^*, 
	\bar{\mbox{\boldmath $ z$}}
	\bigr \rangle_{\mbox{\boldmath \scriptsize $ V$}^*,\mbox{\boldmath \scriptsize $ V$}} 
	& \le & \bigl (\nabla z,\nabla \bar{z} \bigr )_{H^{{\takeshi d}}}+ 
	\bigl (\beta_\varepsilon (z),\bar{z} \bigr )_{H}+
	\varepsilon \bigl (z, \bar{z} \bigr )_{H}+
	\bigl (\nabla _{\Gamma } z_{\Gamma },\nabla _{\Gamma } \bar{z}_{\Gamma } \bigr )_{H_{\Gamma}^{{\takeshi d}} }
	\nonumber \\
	& ~ & 
	+
	\bigl (\beta_{\Gamma ,\varepsilon }(z_{\Gamma }),\bar{z}_{\Gamma } \bigr )_{H_{\Gamma }}
	+\varepsilon \bigl (z_\Gamma , \bar{z}_\Gamma \bigr )_{H_\Gamma }
	\quad \mbox{for all } \bar{\mbox{\boldmath $ z$}}:=(\bar{z},\bar{z}_{\Gamma }) \in \mbox{\boldmath $ V$}.
\end{eqnarray*}
The opposite inequality can be shown as well, by taking $-\delta $ in place of $\delta$. 
Thus, $\partial_* \varphi _\varepsilon $ is single-valued and 
the characterization \eqref{(20)} of $\partial_* \varphi _\varepsilon $ follows. Finally, we see that 
\begin{eqnarray*} 
	&&\int_\Omega \bigl |\beta_\varepsilon (z) \bigr |^2 dx 
	 =  \int_\Omega \frac{1}{\varepsilon ^2} {\revis \bigl |} z -J_\varepsilon (z) {\revis \bigr|}^2 dx \\
	&& \le  \int_\Omega \frac{2}{\varepsilon } \left( 
	  \frac{1}{2\varepsilon } {\revis \bigl |} z-J_\varepsilon (z){\revis \bigr|}^2
	+ \widehat{\beta} {\revis \bigl(} J_\varepsilon (z ){\revis \bigr )} \right) dx
	 = \int_\Omega \frac{2}{\varepsilon } \, \widehat{\beta}_\varepsilon (z)dx 
\end{eqnarray*} 
and 
$$
	\int_\Gamma\bigl |\beta_{\Gamma, \varepsilon}(z_{\Gamma } ) \bigr |^2 d\Gamma
	\le \int_\Gamma \frac{2}{{\takeshi \varrho}  \varepsilon } \, \widehat{\beta}_{\Gamma, \varepsilon} (z_{\Gamma })d\Gamma.
$$
Therefore, the boundedness property in \eqref{(21)} is also true. \hfill $\Box$

\paragraph{Remark 3.1.} 
The estimate \eqref{(21)} is somehow important in order to apply the abstract result in \cite{FuKe}. That was a reason for us to introduce the Moreau-Yosida regularizations $\widehat{\beta}_{\varepsilon }$ and $\widehat{\beta}_{\Gamma,\varepsilon }$, otherwise 
with $\widehat{\beta}$ and $ \widehat{\beta}_{\Gamma }$ instead \eqref{(21)} may not hold.


\bigskip
Now, we recall the fact that 
$$
	\overline{\mbox{\boldmath $ K$}}=
	\mbox{\boldmath $ K$}_{\mbox{\boldmath \scriptsize $ H$}}:=
	\bigl\{ \mbox{\boldmath $ z$} \in \mbox{\boldmath $ H $}\, : \  
	k_* \le (\mbox{\boldmath $ w$},\mbox{\boldmath $ z$})_{\mbox{\boldmath \scriptsize $ H$}}
	\le k^* \bigr\},
$$
is a closed convex subset of $\mbox{\boldmath $ H$}$. 
Then, the following result holds. 

\paragraph{Proposition 3.1.}
{\it For each $\varepsilon \in (0,1]$, there 
exist a unique 
$$\mbox{\boldmath $ u$}_\varepsilon 
\in H^1(0,T;\mbox{\boldmath $ H$}) \cap 
L^\infty (0,T;\mbox{\boldmath $ V$})$$ 
and a pair of functions
$\mbox{\boldmath $ u$}_\varepsilon^* \in L^2(0,T;\mbox{\boldmath $ H$})$ and 
$\lambda _\varepsilon \in L^2(0,T)$ such that 
$$\mbox{\boldmath $ u$}_\varepsilon (t) \in \mbox{\boldmath $ K$}_{\mbox{\boldmath \scriptsize $ H$}} \quad \hbox{for all } t \in [0,T]$$ 
and 
\begin{gather} 
	\mbox{\boldmath $ u$}_\varepsilon '(t)+ \mbox{\boldmath $ u$}_\varepsilon ^* (t) 
	+ \lambda _ \varepsilon (t)\mbox{\boldmath $ w$}
	+\mbox{\boldmath $ \pi $}\bigl( \mbox{\boldmath $ u$}_\varepsilon (t) \bigr)
	= \mbox{\boldmath $ f$}(t)\quad \mbox{in } \mbox{\boldmath $ H$},
	\ \mbox{for a.a.\ } t \in (0,T),  \label{(23)}
\\
	\mbox{\boldmath $ u$}_\varepsilon ^* (t) \in \partial \varphi _\varepsilon 
	\bigl(\mbox{\boldmath $ u$}_\varepsilon (t) \bigr) \quad \mbox{in } 
	\mbox{\boldmath $ H$},
	\ \mbox{for a.a.\ } t \in (0,T),   \label{(24)}
\\
	\lambda _\varepsilon (t) \mbox{\boldmath $ w$} \in 
	\partial I_{\mbox{\boldmath \scriptsize $ K$}_{\mbox{\boldmath \tiny $ H$}}}
	\bigl(\mbox{\boldmath $ u$}_\varepsilon (t) \bigr) \quad \mbox{in } 
	\mbox{\boldmath $ H$},
	\ \mbox{for a.a.\ } t \in (0,T), 
	 \label{(25)}
\\ 
	\mbox{\boldmath $ u$}_\varepsilon (0)=\mbox{\boldmath $ u$}_0 
	\quad \mbox{in } \mbox{\boldmath $ H$}.  \label{(26)}
\end{gather} 
Moreover, $\lambda _\varepsilon$ is given by 
\begin{equation} 
	\lambda _\varepsilon (t) 
	:=\Bigl( 
	\mbox{\boldmath $ f$}(t)-\mbox{\boldmath $ u$}_\varepsilon '(t)-\mbox{\boldmath $ \pi $}
	\bigl(
	\mbox{\boldmath $ u$}_\varepsilon (t) \bigr), \mbox{\boldmath $ z$}_c 
	\Bigr)_{\!\mbox{\boldmath \scriptsize $ H$}}
	-
	\bigl (
	\mbox{\boldmath $ u$}_\varepsilon ^*(t),\mbox{\boldmath $ z$}_c
	\bigr )_{\mbox{\boldmath \scriptsize $ H$}}
	\quad \mbox{for a.a.\ } t \in (0,T),  \label{(27)}
\end{equation} 
where $\mbox{\boldmath $ z$}_c:=(1/{\takeshi \sigma_0} ,1/{\takeshi \sigma_0} ) \in \mbox{\boldmath $ V$}$. 
}%

\paragraph{Proof.} We sketch the basic steps. 
\subparagraph{1.} For a given $\bar{\mbox{\boldmath $ u$}} \in C([0,T];\mbox{\boldmath $ H$})$, there is a unique function $\mbox{\boldmath $ u$} \in H^1(0,T;\mbox{\boldmath $ H$}) \cap L^\infty (0,T;\mbox{\boldmath $ V$})$ solving 
$$
	\mbox{\boldmath $ u$} '(t)+\partial (\varphi _\varepsilon +I_{\mbox{\boldmath \scriptsize $ K$}})
	\bigl( \mbox{\boldmath $ u$} (t) \bigr) 
	\ni \mbox{\boldmath $ f$}(t)
	-\mbox{\boldmath $ \pi $}\bigl( \bar{\mbox{\boldmath $ u$} }(t) \bigr)
	\quad \mbox{in } \mbox{\boldmath $ H$,} 
	\ \, \mbox{for a.a.\ } t \in (0,T),
$$
$$
	\mbox{\boldmath $ u$}(0)=\mbox{\boldmath $ u$}_0 
	\quad \mbox{in } \mbox{\boldmath $ H$}.
$$
Indeed, recalling that  $\mbox{\boldmath $ f$}
	-\mbox{\boldmath $ \pi $}\bigl( \bar{\mbox{\boldmath $ u$} } \bigr)
 \in L^2 (0,T;\mbox{\boldmath $ H$}) $  and $\mbox{\boldmath $ u$}_0 \in D( \varphi _\varepsilon +I_{\mbox{\boldmath \scriptsize $ K$}}) $ (cf.~\eqref{inidata} and \eqref{p3}), 
it suffices to apply, e.g., \cite[Thm.~3.6, p.~72]{Br} for the existence, uniqueness and regularity of the solution $\mbox{\boldmath $ u$} $. Thus, we construct the map 
$$
	\Psi : \bar{\mbox{\boldmath $ u$}} \mapsto \mbox{\boldmath $ u$},
$$
from $C([0,T];\mbox{\boldmath $ H$})$ into itself. 

\subparagraph{2.} For a given pair  
$\bar{\mbox{\boldmath $ u$}}^{(1)}, \bar{\mbox{\boldmath $ u$}}^{(2)} \in C([0,T];\mbox{\boldmath $ H$})$, we can use the estimate (see, e.g., \cite[Lemme~3.1, p.~64]{Br})  
$$
	\bigl| \mbox{\boldmath $ u$}^{(1)}(t) 
	- \mbox{\boldmath $ u$}^{(2)}(t) \bigr|_{\mbox{\boldmath \scriptsize $ H$}}^2
	\le C_{\mbox{\boldmath $ \pi $}}
	\int_{0}^{t} 
	\bigl| \bar{\mbox{\boldmath $ u$}}^{(1)}(s ) 
	-\bar{\mbox{\boldmath $ u$}}^{(2)}(s ) \bigr|_{\mbox{\boldmath \scriptsize $ H$}}^2 ds 
	\quad \mbox{ for all } t \in [0,T],
$$ 
where  ${\mbox{\boldmath $ u$}}^{(i)} = \Psi (\bar{\mbox{\boldmath $ u$}}^{(i)})$, $i=1,2$,
and $C_{\mbox{\boldmath $ \pi $}}>0$ is a  positive constant depending only  on $L$ and $L_{\Gamma}$ (cf.~\eqref{pilip}). Then, by recurrence one shows that there exists a suitable $k \in \mathbb{N}$ such that $\Psi^k$ is a 
contraction mapping in $C([0,T];\mbox{\boldmath $ H$})$, and consequently there 
exists a unique solution $\mbox{\boldmath $ u$}_\varepsilon $ of the 
problem \eqref{(18)}--\eqref{(19)}. 

\subparagraph{3.} Now, in order to conclude the proof we can just apply Theorem~2.3 and Remark~3 of \cite{FuKe}. In fact,  
in view of Lemma 3.1, it is not difficult to check the validity 
of the assumptions (A1)--(A5) of \cite{FuKe} in our case. In particular, let us point out that the coercivity property stated in \cite[(A5)]{FuKe} comes from the definition of $\varphi_\varepsilon$. However, one important point regards the density of 
$\mbox{\boldmath $ V$}$ in $\mbox{\boldmath $ H$}$, for which we refer the reader to the  
Appendix. Finally, we use the fact that 
$\overline{\mbox{\boldmath $ K$}}= \mbox{\boldmath $ K$}_{\mbox{\boldmath \scriptsize $ H$}}$. 
\hfill $\Box$

\bigskip
Thanks to Proposition~3.1 and Lemma~3.1, 
we arrive at the following weak formulation of~\eqref{(23)}:
\begin{align} 
	\lefteqn{ 
	\int_{\Omega}^{} g_\varepsilon  (t)z dx 
	+ \int_{\Gamma }^{} g_{\Gamma,\varepsilon } (t) z_\Gamma d\Gamma
	} \nonumber \\
	& = 
	\int_{\Omega }^{} \nabla u_\varepsilon (t) \cdot \nabla z dx 
	+ \int_{\Omega }^{} \beta_\varepsilon \bigl( u_\varepsilon (t) \bigr)z dx 
	+ \varepsilon \int_{\Omega }^{}  u_\varepsilon (t) 
	z dx 
	+ \int_{\Gamma }^{} \nabla _\Gamma u_{\Gamma, \varepsilon }(t)\cdot 
	\nabla _\Gamma z_\Gamma d\Gamma \nonumber \\
	& \qquad{}
	+ \int_{\Gamma }^{} \beta_{\Gamma, \varepsilon} \bigl( u_{\Gamma, \varepsilon}(t) \bigr)
	z_\Gamma d\Gamma 
	+  \varepsilon \int_{\Gamma }^{} u_{\Gamma, \varepsilon}(t)
	z_\Gamma d\Gamma 
	\quad \mbox{for all } \mbox{\boldmath $ z$}:=(z, z_\Gamma ) \in \mbox{\boldmath $ V$},
	 \label{(28)}
\end{align} 
where  
$$
	g_\varepsilon :=f_\varepsilon 
	-u_\varepsilon '
	-\lambda _\varepsilon w 
	-\pi (u_\varepsilon )
	\in L^2(0,T;H),
$$
$$
	g_{\Gamma, \varepsilon} :=f_{\Gamma, \varepsilon} 
	-u_{\Gamma, \varepsilon}'
	-\lambda _\varepsilon w _{\Gamma }
	-\pi_\Gamma  (u_{\Gamma, \varepsilon} )
	\in L^2(0,T;H_\Gamma ).
$$
Moreover, we point out the following regularity properties for the solution. 

\paragraph{Proposition 3.2.} {\it 
For each $\varepsilon \in (0,1]$ we have that
$u_\varepsilon \in L^2 (0,T;H^2 (\Omega ))$ and 
$u_{\Gamma, \varepsilon} \in L^2 ( 0,T;H^2(\Gamma))$. }

\paragraph{Proof.}
First, we take $z \in {\mathcal D}(\Omega )$, which entails that $z_\Gamma =0$, in \eqref{(28)} and get 
$$
	-\Delta u_\varepsilon (t) 
	=g_\varepsilon (t) -
	\beta_\varepsilon \bigl( u_\varepsilon (t) \bigr) -\varepsilon u_\varepsilon (t) 
	\quad \mbox{in } {\mathcal D}'(\Omega), \quad \mbox{for a.a.\ } t \in (0,T).
$$
This implies that $-\Delta u_\varepsilon \in L^2(0,T;H)$ due to the regularity of the right hand side. On the other hand, we already know that 
$u_\varepsilon \in L^\infty (0,T;V)$ and $u_{\Gamma, \varepsilon } \in L^\infty (0,T;V_\Gamma )$. Then, we infer that (see, e.g.,\ \cite[Thm.~3.2, p.~1.79]{BrGi})
$$
	u_\varepsilon \in L^2 \bigl (0,T;H^{3/2} (\Omega ) \bigr )
$$
and consequently, by a trace theorem \cite[Thm.~2.27, p.~1.64]{BrGi}, $\partial_{\nu} u_\varepsilon \in L^2(0,T;H_\Gamma )$.
At this point, from the variational equality \eqref{(28)} we can obtain 
the characterization on the boundary
$$
	- \Delta_\Gamma u_{\Gamma, \varepsilon } 
	= g_{\Gamma, \varepsilon}
	- \partial_{\nu } u_\varepsilon 
	- \beta_{\Gamma, \varepsilon} (u_{\Gamma, \varepsilon} )
	- \varepsilon u_{\Gamma, \varepsilon} 
	\quad \mbox{a.e.\ on } \Sigma, 
$$
and the information that $\Delta_{\Gamma } u_{\Gamma, \varepsilon} \in L^2(0,T;H_\Gamma )$ 
implies (see, e.g., \cite[p.~104]{Gr})
$$
	u_{\Gamma, \varepsilon} \in L^2 \bigl ( 0,T;H^2(\Gamma) \bigr ).
$$
Finally, this yields in particular that $u_{\Gamma, \varepsilon} \in L^2 ( 0,T;H^{3/2}(\Gamma))$, whence (quoting again \cite[Thm.~3.2, p.~1.79]{BrGi})
$$
	u_\varepsilon \in L^2 \bigl (0,T;H^2 (\Omega ) \bigr ). 
$$
\vskip-0.5cm
\hfill $\Box$ 

\bigskip
By virtue of this lemma, our approximate problem can be written as 
\begin{gather} 
	\frac{\partial u_\varepsilon }{\partial t}-\Delta u_\varepsilon 
	+ \beta_\varepsilon(u_\varepsilon) 
	+ \varepsilon u_\varepsilon 
	+ \pi(u_\varepsilon) 
	+ \lambda_\varepsilon w = f
	\quad \mbox{a.e.\ in } Q,  \label{(29)}
\\ 
	\partial_\nu u_\varepsilon
	+ \frac{\partial u_{\Gamma, \varepsilon  }}{\partial t}
	- \Delta_{\Gamma } u_{\Gamma, \varepsilon }
	+ \beta_{\Gamma, \varepsilon} (u_{\Gamma, \varepsilon })
	+ \varepsilon u_{\Gamma, \varepsilon }
	+ \pi_{\Gamma }(u_{\Gamma, \varepsilon}) 
	+ \lambda_\varepsilon w _{\Gamma } = f_\Gamma 
	\quad \mbox{a.e.\ on } \Sigma,  \label{(30)}
\\ 
	u_{\Gamma, \varepsilon } = {u_\varepsilon}_{| _\Gamma} 
	\quad \mbox{a.e.\ on } \Sigma,  \label{(31)}
\\
	u_\varepsilon (0)=u_0 \quad \mbox{a.e.\ in } \Omega, \quad 
	u_{\Gamma, \varepsilon}(0)=u_{0\Gamma } \quad \mbox{a.e.\ on } \Gamma,
	 \label{(32)}
\\ 
	k_* \le k_\varepsilon (t):=\int_{\Omega }^{}w u_\varepsilon (t) dx 
	+ \int_{\Gamma }^{} w_{\Gamma } u_{\Gamma, \varepsilon}(t) d\Gamma 
	\le k^* 
	\quad \mbox{for all }t \in [0,T],  \label{(33)}
\\
	\lambda _\varepsilon(t) \in \partial I_{[k_*,k^*]}{\fukao \bigl( }k_\varepsilon (t) {\fukao \bigr)}
	\quad \mbox{for a.a.\ } t \in (0,T).  \label{(34)}
\end{gather}

\paragraph{Remark 3.2.} As 
$\mbox{\boldmath $ u$}_\varepsilon (t) \in \mbox{\boldmath $ K$}_{\mbox{\boldmath \scriptsize $ H$}}$ for all $t\in [0,T]$, 
we claim that the last condition is equivalent to \eqref{(25)}. Actually, let us assume \eqref{(34)}. For each $\mbox{\boldmath $ z$} \in \mbox{\boldmath $ K$}_{\mbox{\boldmath \scriptsize  $ H$}}$, 
there exist uniquely $\alpha \in \mathbb{R}$ with 
$k_* \le \alpha \le k^*$ and $\mbox{\boldmath $ z$}_N \in \mbox{\boldmath $ H$}$ with 
$(\mbox{\boldmath $ w$},\mbox{\boldmath $ z$}_N)_{\mbox{\boldmath \scriptsize $ H$}}=0$ such that 
$$
	\mbox{\boldmath $ z$}=\alpha \mbox{\boldmath $ z$}_c+\mbox{\boldmath $ z$}_N, 
	\quad (\mbox{\boldmath $ w$},\mbox{\boldmath $ z$})_{\mbox{\boldmath \scriptsize $ H$}}=\alpha.
$$
Therefore, from the definition of subdifferential it follows that 
$\lambda _\varepsilon (t) \cdot (  
(\mbox{\boldmath $ w$},\mbox{\boldmath $ z$})_{\mbox{\boldmath \scriptsize $ H$}}
- k_\varepsilon (t)) \le 0$, 
namely 
\begin{equation}
	\bigl( 
	\lambda _\varepsilon (t)\mbox{\boldmath $ w$},
	\mbox{\boldmath $ z$}-
	\mbox{\boldmath $ u$}_\varepsilon (t) 
	\bigr )_{\mbox{\boldmath \scriptsize $ H$}} \le 0
	\quad \mbox{for all } \mbox{\boldmath $ z$} \in \mbox{\boldmath $ K$}_{\mbox{\boldmath \scriptsize $ H$}}. \label{p5} 
\end{equation}
Thus, \eqref{(25)} holds. On the other hand, let us assume \eqref{(25)}. Then, we can take  
 $\mbox{\boldmath $ z$}:=r\mbox{\boldmath $ z$}_c$,  $r \in [k_*,k^*]$, as test function
in \eqref{p5}, so that we obtain \eqref{(34)} by recalling that  $
	k_\varepsilon (t) = (\mbox{\boldmath $ w$},
	\mbox{\boldmath $ u$}_\varepsilon (t) 
	)_{\mbox{\boldmath \scriptsize $ H$}}$ from \eqref{(33)}.

\subsection{A priori estimates}

In this subsection, we obtain the uniform estimates independent of $\varepsilon >0$. 
\paragraph{Lemma 3.2.}
{\it There exist a positive constant $M_1$, 
independent of $\varepsilon \in (0,1]$, such that 
\begin{align}
	&|u_\varepsilon|_{H^1(0,T;H)}+|u_\varepsilon |_{L^\infty (0,T;V)}
	+ \sup_{t \in (0,T)} \int_{\Omega }^{} \widehat{\beta}_\varepsilon 
	{\revis \bigl (}u_\varepsilon (t) {\revis \bigr) }   
	dx\nonumber \\
	&{}	+
	|u_{\Gamma, \varepsilon}|_{H^1(0,T;H_\Gamma )}+|u_{\Gamma, \varepsilon}|_{L^\infty (0,T;V_\Gamma )}
	+ \sup_{t \in (0,T)} \int_{\Gamma}^{} \widehat{\beta}_{\Gamma, \varepsilon} 
	{\revis \bigl (}u_{\Gamma, \varepsilon} (t){\revis \bigr) } d\Gamma  
	\le M_1. \label{p6}
\end{align} 
}%

\paragraph{Proof.} We can add $u_\varepsilon $ to 
both sides of \eqref{(29)} and  $u_{\Gamma ,\varepsilon }$ to 
both sides of \eqref{(30)}, then test \eqref{(29)} by 
${\pier (\partial u_\varepsilon/\partial t) }\in L^2(0,T;H)$ and 
use boundary conditions \eqref{(30)}--\eqref{(31)}. {\pier Then, we deduce that}
\begin{eqnarray*}
	\lefteqn{ 
	\int_{0}^{t} \bigl |u_\varepsilon '(s ) \bigr |_H^2 ds 
	+ \frac{1}{2} \bigl |u_\varepsilon (t) \bigr |_V^2 
	+ \int_{\Omega }^{} \widehat{\beta}_\varepsilon 
	\bigl (u_\varepsilon (t) \bigr) dx
	+ \frac{\varepsilon }{2} \bigl |u_\varepsilon (t) \bigr |_H^2
	} \nonumber \\
	& ~ & {} 
	+ \int_{0}^{t} \bigl |u_{\Gamma, \varepsilon }'(s ) \bigr |_{H_\Gamma }^2 ds 
	+ \frac{1}{2} \bigl |u_{\Gamma, \varepsilon }(t) \bigr |_{V_\Gamma }^2 
	+ \int_{\Gamma }^{} \widehat{\beta}_{\Gamma, \varepsilon} 
	\bigl (u_{\Gamma, \varepsilon} (t) \bigr) d\Gamma
	+ \frac{\varepsilon }{2} \bigl |u_{\Gamma, \varepsilon }(t) \bigr |_{H_\Gamma}^2 
	\nonumber \\
	& ~ & {}
	+ \int_{0}^{t} \lambda _\varepsilon (s ) 
	\left\{ 
	\int_{\Omega }^{}
	w u_\varepsilon '(s)dx 
	+ \int_{\Gamma }^{}
	w_\Gamma u_{\Gamma, \varepsilon} '(s )d\Gamma 
	\right\} ds \nonumber \\
	& \le & 
	\frac{1}{2} |u_0|_V^2 + 
	\int_{\Omega }^{}
	\widehat{\beta}_\varepsilon (u_0) dx  
	+\frac{\varepsilon }{2} |u_0|_H^2
	+ \frac{1}{2} |u_{0\Gamma}|_{V_\Gamma }^2 
	+ \int_{\Gamma }^{} \widehat{\beta}_{\Gamma,\varepsilon  }
	(u_{0\Gamma} ) d\Gamma 
	+ \frac{\varepsilon }{2} |u_{0\Gamma} |_{H_\Gamma }^2 
	\nonumber \\
	& ~ & {}
	+ \int_{0}^{t} \Bigl ( f(s )-\pi 
	\bigl (u_\varepsilon (s ) 
	\bigr ) ,u_\varepsilon '(s ) 
	\Bigr )_{\!H} ds 
	+
	\int_{0}^{t} \Bigl ( f_\Gamma (s )-\pi_\Gamma 
	\bigl (u_{\Gamma , \varepsilon} (s ) 
	\bigr ) ,u_{\Gamma, \varepsilon} '(s ) 
	\Bigr )_{\!H_\Gamma } ds, 
\end{eqnarray*} 
for all $t \in [0,T]$. 
We {\pier point out} that~{\pier (cf.~\eqref{p4})}
$$
	\int_{\Omega }^{}
	\widehat{\beta}_\varepsilon (u_0) dx 
	\le 
	\int_{\Omega }^{}
	\widehat{\beta}(u_0) dx<+\infty , \quad 
	\int_{\Gamma }^{} \widehat{\beta}_{\Gamma,\varepsilon  }
	(u_{0\Gamma} ) d\Gamma 
	\le
	\int_{\Gamma }^{} \widehat{\beta}_\Gamma
	(u_{0\Gamma} ) d\Gamma <+ \infty.
$$
Also note that, from \eqref{(33)}, \eqref{(34)} and 
the chain rule differentiation lemma (see, {\pier e.g.,
\cite[Lemma~4.4, p.~158]{Ba} or \cite[Lemme~3.3, p.~73]{Br}}), 
the last term on the left hand side is exactly 
$$
	\int_{0}^{t} \lambda _\varepsilon (s) k_\varepsilon '(s) ds 
	= I_{[k_*,k^*]} {\fukao \bigl(}k_\varepsilon(t) {\fukao \bigr)} -I_{[k_*,k^*]}(k_0) \equiv 0 
	\quad \mbox{for all }  t \in [0,T],
$$
where $k_0:=(\mbox{\boldmath $ w$},\mbox{\boldmath $ u$}_0)_{\mbox{\boldmath \scriptsize $ H$}}$.
Moreover, 
there exists a positive constant $\tilde{M}_1$, {\pier depending only on 
$L$, $L_\Gamma$, $|u_0|_H$, $|u_{0\Gamma }|_{H_\Gamma}$ and $T$,} such that 
$$
	\int_{0}^{t} \Bigl ( f(s )-\pi 
	\bigl (u_\varepsilon (s ) 
	\bigr ) ,u_\varepsilon '(s ) 
	\Bigr )_{\!H} ds \le \frac{1}{2} \int_{0}^{t} {\fukao \bigl|} u_\varepsilon '(s) {\fukao \bigr|}_H^2 ds 
	+ \tilde{M}_1 \int_{0}^{t} 
	\Bigl( 1 +  \bigl |f(s) \bigr |_H^2 
	+ \bigl | u_\varepsilon (s) \bigr |_H^2 \Bigr) 
	ds,
$$
and 
\begin{eqnarray*}
	\lefteqn{ 
	\int_{0}^{t}  \Bigl ( f_\Gamma (s )-\pi_\Gamma 
	\bigl (u_{\Gamma , \varepsilon} (s ) 
	\bigr ) ,u_{\Gamma, \varepsilon} '(s ) 
	\Bigr )_{\!H_\Gamma } ds 
	} \\
	& \le & \frac{1}{2} \int_{0}^{t} {\fukao \bigl|}u_{\Gamma, \varepsilon} '(s){\fukao \bigr|}_{H_\Gamma}^2 ds 
	+ \tilde{M}_1 \int_{0}^{t} 
	\Bigl( 1 +  \bigl |f_\Gamma (s) \bigr |_{H_\Gamma }^2 
	+ \bigl | u_{\Gamma, \varepsilon }(s) \bigr |_{H_\Gamma}^2 \Bigr) 
	ds  
\end{eqnarray*} 
{\pier for all $t \in [0,T]$.}
Collecting the estimates and applying the Gronwall inequality, 
we {\pier easily get \eqref{p6}}. \hfill $\Box$

\bigskip
Thanks to the growth conditions \eqref{(5)}--\eqref{(6)} {\pier (see also \eqref{(15)}--\eqref{(16)}), we obtain the following bound.} 

\paragraph{Lemma 3.3.}
{\it There exist a positive constant $M_2$, 
independent of $\varepsilon \in (0,1]$, such that} 
$$
	|\lambda _\varepsilon |_{L^2(0,T)} \le M_2.
$$

\paragraph{Proof.} From the expression of $\lambda _\varepsilon $, given at \eqref{(27)}, 
we have that 
\begin{align*} 
	&\lambda _\varepsilon (t)
	 =  \frac{1}{{\takeshi \sigma_0} } \int_{\Omega }^{}
	\Bigl( f(t)-u_\varepsilon '(t)
	-\pi \bigl (u_\varepsilon (t) \bigr )
	\Bigr) dx \\
	&\qquad {}
	+
	\frac{1}{{\takeshi \sigma_0} } \int_{\Gamma }^{}
	\Bigl( f_\Gamma (t)-u_{\Gamma, \varepsilon} '(t)
	-\pi_\Gamma  \bigl (u_{\Gamma, \varepsilon} (t) \bigr )
	\Bigr) d\Gamma 
	- \bigl( \mbox{\boldmath $ u$}^*_\varepsilon (t),
	\mbox{\boldmath $ z$}_c \bigr)_{\mbox{\boldmath \scriptsize $ H$}} 
	\quad \mbox{for a.a.\ } t \in (0,T).
\end{align*}
{\pier As $\mbox{\boldmath $ z$}_c  =(1/{\takeshi \sigma_0} ,1/{\takeshi \sigma_0} ) \in \mbox{\boldmath $ V$}$ 
and $\mbox{\boldmath $ u$}_\varepsilon ^*(t) 
= \partial \varphi _\varepsilon (\mbox{\boldmath $ u$}_\varepsilon (t))$ in 
$\mbox{\boldmath $ H$}$, using \eqref{(20)} we obtain}
$$
	\bigl( \mbox{\boldmath $ u$}^*_\varepsilon (t),
	\mbox{\boldmath $ z$}_c \bigr)_{\mbox{\boldmath \scriptsize $ H$}} 
	= \frac{1}{{\takeshi \sigma_0} } \int_{\Omega }^{} \Bigl( 
	\beta_\varepsilon \bigl( u_\varepsilon (t) \bigr ) 
	+ \varepsilon u_\varepsilon (t) \Bigr) dx 
	+ \frac{1}{{\takeshi \sigma_0} } \int_{\Gamma }^{} 
	\Bigl( \beta_{\Gamma, \varepsilon }\bigl( u_{\Gamma,\varepsilon }(t) \bigr ) 
	+ \varepsilon u_{\Gamma, \varepsilon }(t)
	\Bigr) d\Gamma
$$
{\pier for a.a. $t\in (0,T)$.} Then, we can estimate $\lambda _\varepsilon $ as follows:
\begin{eqnarray*} 
	\lefteqn{ 
	|\lambda _\varepsilon|_{L^2(0,T)}^2
	} \nonumber \\
	& \le & \tilde{M}_2 
	\Bigl(1+ |f|_{L^2(0,T;H)}^2 +|u_\varepsilon|_{H^1(0,T;H)}^2
	\Bigr) 
	+ \tilde{M}_2 
	\Bigl(1+ |f_\Gamma |_{L^2(0,T;H_\Gamma )}^2 +|u_{\Gamma , \varepsilon}|_{H^1(0,T;H_\Gamma )}^2
	\Bigr) \\
	& ~ & 
	{}+ \tilde{M}_2 
	\sup_{t \in (0,T)} 
	\left( 
	\left| 
	\int_{\Omega }^{} 
	\beta_\varepsilon \bigl( u_\varepsilon (t) \bigr) dx 
	\right|^2 
	+
	\left| 
	\int_{\Gamma }^{} 
	\beta_{\Gamma, \varepsilon} \bigl( u_{\Gamma, \varepsilon} (t) \bigr) d\Gamma 
	\right|^2 
	\right),
\end{eqnarray*}
where $\tilde{M}_2$ is a positive constant, {\pier depending on 
${\takeshi \sigma_0} $, $|\Omega|$, $|\Gamma|$, $L$, $L_\Gamma$, $|u_0|_H$, 
$|u_{0\Gamma}|_{H_\Gamma }$ and $T$}. 
Now, we use {\pier the properties \eqref{(15)}--\eqref{(16)} along with the estimate \eqref{p6}} to conclude. \hfill $\Box$

\paragraph{Lemma 3.4.}
{\it There exist {\pier two} positive constants $M_3$ and $M_4$, 
independent of $\varepsilon \in (0,1]$, such that 
$$
	\bigl| \beta_\varepsilon (u_\varepsilon ) \bigr |_{L^2(0,T;H)} 
	+ \bigl |\beta_\varepsilon (u_{\Gamma,\varepsilon }) \bigr |_{L^2(0,T;H_\Gamma )} \le M_3,
$$
$$
	|u_\varepsilon |_{L^2(0,T;H^{3/2}(\Omega ))}
	+ 
	|\partial_\nu u_\varepsilon |_{L^2(0,T;H_\Gamma )} 
	\le M_4.
$$
}%

\paragraph{Proof.} Testing \eqref{(29)} by 
$\beta_\varepsilon (u_\varepsilon ) \in L^2(0,T;{\pier V})$ 
and using \eqref{(30)}--\eqref{(31)}{\pier , we infer that}
\begin{eqnarray*}
	\lefteqn{ 
	\int_{\Omega }^{} \widehat{\beta}_\varepsilon 
	\bigl( u_\varepsilon (t) \bigr) dx 
	+ \int_0^t\!\!\! \int_{\Omega }^{} \beta '_\varepsilon 
	\bigl( u_\varepsilon (s) \bigr)
	\bigl | \nabla u_\varepsilon (s) \bigr |^2 dxds
	+ 
	\int_{0}^{t} 
	\Bigl | \beta_\varepsilon 
	\bigl (u_\varepsilon (s) \bigr) \Bigr |^2_H ds 
	} \nonumber \\
	& ~ & {} 
	+ \varepsilon \int_0^t\!\!\! 
	\int_{\Omega }^{} u_\varepsilon (s) \beta_\varepsilon 
	\bigl( u_\varepsilon (s) \bigr) dx ds 
	+ \int_{\Gamma }^{} \widehat{\beta}_\varepsilon 
	\bigl( u_{\Gamma ,\varepsilon} (t) \bigr) d\Gamma 
	\nonumber \\
	& ~ & {}
	+ \int_0^t\!\!\! \int_{\Gamma }^{} \beta '_\varepsilon 
	\bigl( u_{\Gamma, \varepsilon} (s) \bigr)
	\bigl | \nabla_\Gamma  u_{\Gamma,\varepsilon} (s) \bigr |^2 d\Gamma ds
	+ 
	\int_0^t\!\!\! \int_{\Gamma }^{} 
	\beta_\varepsilon \big( u_{\Gamma, \varepsilon } (s)\bigr)
	\beta_{\Gamma, \varepsilon } \bigl( u_{\Gamma, \varepsilon }(s) \bigr)
	d\Gamma ds
	\nonumber \\
	& ~ & {} 
	+ \varepsilon \int_0^t\!\!\! 
	\int_{\Gamma }^{} 
	u_{\Gamma, \varepsilon } (s)
	\beta_\varepsilon \bigl (u_{\Gamma, \varepsilon} (s)\bigr )
	d\Gamma ds \\
	& \le & 
	\int_{\Omega }^{}
	\widehat{\beta}_\varepsilon (u_0) dx  
	+ \int_{\Gamma }^{} \widehat{\beta}_\varepsilon (u_{0\Gamma }) d\Gamma 
	+ \int_{0}^{t} \Bigl ( f(s )-\pi 
	\bigl (u_\varepsilon (s ) 
	\bigr ) - \lambda _\varepsilon w ,\beta_\varepsilon 
	\bigl( u_\varepsilon (s ) \bigr) 
	\Bigr )_{\!H} ds 
	\\
	& ~ & {}
	+
	\int_{0}^{t} \Bigl ( f_\Gamma (s )-\pi_\Gamma 
	\bigl (u_{\Gamma , \varepsilon} (s ) 
	\bigr ) 
	-\lambda _\varepsilon w_\Gamma ,
	\beta_\varepsilon 
	\bigl( u_{\Gamma, \varepsilon}(s ) 
	\bigr) 
	\Bigr )_{\!H_\Gamma } ds
\end{eqnarray*} 
for all $t \in [0,T]$. Now, we use the {\pier property}~\eqref{(17)} to deduce that 
\begin{eqnarray*} 
\lefteqn{ 
	\int_0^t\!\!\! \int_{\Gamma }^{} 
	\beta_\varepsilon \big( u_{\Gamma, \varepsilon } (s)\bigr)
	\beta_{\Gamma, \varepsilon } \bigl( u_{\Gamma, \varepsilon }(s) \bigr)
	d\Gamma ds 
	} \nonumber \\
	& = & 
	\int_0^t\!\!\! \int_{\Gamma }^{} 
	\Bigl| \beta_\varepsilon \big( u_{\Gamma, \varepsilon } (s)\bigr) 
	\Bigr| 
	\Bigl|
	\beta_{\Gamma, \varepsilon } \bigl( u_{\Gamma, \varepsilon }(s) \bigr)
	\Bigr| 
	d\Gamma ds 
	\\
	& \ge & 
	\frac{1}{{\takeshi \varrho}  } \int_0^t\!\!\! \int_{\Gamma }^{} 
	\Bigl | \beta_\varepsilon \bigr( u_{\Gamma ,\varepsilon } (s) \bigr) \Bigr | ^2 
	d\Gamma ds
	- \frac{c_0}{{\takeshi \varrho}  }\int_0^t\!\!\! \int_{\Gamma }^{} 
	\Bigl |\beta_\varepsilon \bigr( u_{\Gamma ,\varepsilon } (s) \bigr) \Bigr| 
	d\Gamma ds
\end{eqnarray*} 
for all $t \in [0,T]$, 
because $\beta_\varepsilon (r)$ and $\beta_{\Gamma, \varepsilon }(r)$ have the same sign 
for all $r \in \mathbb{R}$. 
We also note that 
$$
	{\pier \int_0^t\!\!\! \int_{\Omega }^{} \beta '_\varepsilon 
	\bigl( u_\varepsilon (s) \bigr)
	\bigl | \nabla u_\varepsilon (s) \bigr |^2 dxds \ge 0,  \ \quad
	\varepsilon \int_0^t\!\!\! \int_{\Omega }^{} u_\varepsilon (s) 
	\beta_\varepsilon \bigl( u_\varepsilon (s) \bigr) dxds \ge 0,}
$$
$$
	{\pier \int_0^t\!\!\! \int_{\Gamma }^{} \beta '_\varepsilon 
	\bigl( u_{\Gamma, \varepsilon} (s) \bigr)
	\bigl | \nabla_\Gamma  u_{\Gamma,\varepsilon} (s) \bigr |^2 d\Gamma ds \ge 0,
	\quad\ 	\varepsilon \int_0^t\!\!\! \int_{\Gamma }^{} u_{\Gamma, \varepsilon} (s) 
	\beta_\varepsilon \bigl( u_{\Gamma, \varepsilon} (s) \bigr) d\Gamma ds \ge 0}
$$
for all $t \in [0,T]$ and{\pier , in view of \eqref{inidata} and \eqref{p4},}
\begin{align*} 
	&\int_{\Gamma }^{} \widehat{\beta}_\varepsilon (u_{0\Gamma }) d\Gamma
		{\pier  {}={}} \int_{\Gamma }^{} 
	\int_{0}^{u_{0\Gamma }} \beta_\varepsilon (r) dr d\Gamma
	 \le  \int_{\Gamma }^{} 
	\int_{0}^{u_{0\Gamma }} {\takeshi \varrho}  \beta_{\Gamma, \varepsilon} (r) dr d\Gamma 
	+ \int_{\Gamma }^{} c_0 |u_{0\Gamma }| d \Gamma \\
	& {}\le {\takeshi \varrho}  \int_{\Gamma }^{} 
	\widehat{\beta} _{\Gamma, \varepsilon} (u_{0\Gamma }) d\Gamma 
	+ c_0 |u_{0\Gamma }|_{L^1(\Gamma )}
	\le {\takeshi \varrho}  \int_{\Gamma }^{} 
	\widehat{\beta} _{\Gamma} (u_{0\Gamma }) d\Gamma 
	+ c_0 |u_{0\Gamma }|_{L^1(\Gamma )}
	 <  +\infty.
\end{align*} 
Moreover, there exists a positive constant $\tilde{M}_3$, 
{\pier depending} on ${\takeshi \varrho} $, $L$, $L_\Gamma$, $|u_0|_H$, $|u_{0\Gamma}|_{H_\Gamma }
${\pier, $T$ and} independent of $\varepsilon \in (0,1]$, such that
\begin{eqnarray*} 
	\lefteqn{ 
	\int_{0}^{t} \Bigl ( f(s )-\pi 
	\bigl (u_\varepsilon (s ) 
	\bigr ) - \lambda _\varepsilon(s) w ,\beta_\varepsilon 
	\bigl( u_\varepsilon (s ) \bigr) 
	\Bigr )_{\!H} ds 
	} \nonumber \\
	& \le & 
	\frac{1}{2} \int_{0}^{t} 
	\Bigl| \beta_\varepsilon 
	\bigl( u_\varepsilon (s ) \bigr) 
	\Bigr|_{\!H}^2 ds
	+ 
	\tilde{M}_3
	\left( 
	1+|f|_{L^2(0,T;H)}^2 + |u_\varepsilon |_{L^2(0,T;H)}^2 
	 + |\lambda _\varepsilon |_{L^2(0,T)}^2 |w|_H^2
	\right),
	\\
	\lefteqn{ 
	\int_{0}^{t} \Bigl ( f_\Gamma (s )-\pi_\Gamma 
	\bigl (u_{\Gamma , \varepsilon} (s ) 
	\bigr ) 
	-\lambda _\varepsilon(s) w_\Gamma ,
	\beta_\varepsilon 
	\bigl( u_{\Gamma, \varepsilon}(s ) 
	\bigr) 
	\Bigr )_{\! H_\Gamma } ds
	} \nonumber \\
	& \le & 
	\frac{1}{2{\takeshi \varrho}  } \int_{0}^{t} 
	\Bigl| \beta_\varepsilon 
	\bigl( u_{\Gamma, \varepsilon }(s ) \bigr) 
	\Bigr|_{H_\Gamma }^2 ds  \\
	&&+ 
	\tilde{M}_3
	\left( 
	1+|f_\Gamma |_{L^2(0,T;H_\Gamma )}^2 + |u_{\Gamma,\varepsilon} |_{L^2(0,T;H)}^2 
	 + |\lambda _\varepsilon |_{L^2(0,T)}^2 |w_\Gamma |_{H_\Gamma }^2
	\right),
\end{eqnarray*} 
for all $t \in [0,T]$. Thus, we deduce that 
there {\pier is a positive constant $M_3$, which 
depends only} on 
$|f|_{L^2(0,T;H)}$, 
$|f_\Gamma|_{L^2(0,T;H_\Gamma )}$, 
$|u_0|_H$, 
$|u_{0\Gamma}|_{H_\Gamma}$, $M_1$ and $M_2$, 
{\pier such that} 
$$
	\bigl |\beta_\varepsilon (u_\varepsilon )\bigr |_{L^2(0,T;H)} 
	+ \bigl |\beta_\varepsilon (u_{\Gamma,\varepsilon })\bigr |_{L^2(0,T;H_\Gamma )} \le M_3.
$$
Now, we can compare the terms in \eqref{(29)} and 
conclude that 
$$
	|\Delta u_\varepsilon |_{L^2(0,T;H)} \le M_4, 
$$
whence, {\pier recalling \eqref{p6}} and applying the theory of the elliptic regularity (see, e.g.{\pier ,\ \cite[Thm.~3.2, p.~1.79]{BrGi}}){\pier , we have that
}
$$
	|u_\varepsilon |_{L^2(0,T;H^{3/2}(\Omega ))}
	\le M_4
$$
and, {\pier owing to} the trace theory (see, e.g.{\takeshi ,\ \cite[Thm.~2.25, 
p.~1.62]{BrGi}}){\pier , that} 
$$
	|\partial_\nu u_\varepsilon |_{L^2(0,T;H_\Gamma )} 
	\le M_4
$$
{\pier for some constant $M_4$ independent of $\varepsilon \in (0,1]$.}
\hfill $\Box$

\paragraph{Lemma 3.5.}
{\it There exist {\pier three} positive constants $M_5$, $M_6$ and $M_7$,
independent of $\varepsilon \in (0,1]$, such that 
$$
	\bigl |\beta_{\Gamma, \varepsilon} (u_{\Gamma,\varepsilon }) 
	\bigr |_{L^2(0,T;H_\Gamma )} \le M_5,
$$
$$
	|u_{\Gamma, \varepsilon} |_{L^2(0,T;H^2(\Gamma))}
	\le M_6,
$$
$$
	|u_\varepsilon |_{L^2(0,T;H^2(\Omega ))}
	\le M_7.
$$
}%

\paragraph{Proof.} We test \eqref{(30)} by 
$\beta_{\Gamma, \varepsilon} (u_{\Gamma, \varepsilon}) \in L^2(0,T;H_\Gamma)$ and 
integrate on the boundary, {\pier deducing that} 
\begin{eqnarray}
	\lefteqn{ 
	\int_{\Gamma}^{} \widehat{\beta}_{\Gamma, \varepsilon}
	\bigl( u_{\Gamma, \varepsilon} (t) \bigr) d\Gamma 
	+ \int_0^t\!\!\! \int_{\Gamma}^{} \beta '_{\Gamma, \varepsilon }
	\bigl( u_{\Gamma, \varepsilon} (s) \bigr)
	\bigl | \nabla_\Gamma  u_\varepsilon (s) \bigr |^2 d\Gamma ds
	} \nonumber \\
	& ~ & {}
	+ \int_{0}^{t}\Bigl| \beta_{\Gamma, \varepsilon }
	\bigl( u_{\Gamma, \varepsilon}(s) 
	\bigr) \Bigr|_{H_\Gamma }^2  ds
	+ \varepsilon \int_0^t\!\!\! \int_{\Gamma }^{} u_{\Gamma, \varepsilon } (s) 
	\beta_{\Gamma, \varepsilon } \bigl (u_{\Gamma, \varepsilon }(s) \bigr) d\Gamma ds 
	\nonumber\\
	& \le & 
	\int_{\Gamma}^{} \widehat{\beta}_{\Gamma, \varepsilon}
	\bigl( u_{0\Gamma} \bigr) d\Gamma \nonumber \\
	& ~ & {}
	+
	\int_{0}^{t} 
	\Bigl ( f_\Gamma (s )-\partial_\nu u_{\Gamma , \varepsilon}(s)
	-\pi_\Gamma 
	\bigl (u_{\Gamma , \varepsilon} (s ) 
	\bigr ) 
	-\lambda _\varepsilon w_\Gamma ,
	\beta_{\Gamma, \varepsilon }
	\bigl( u_{\Gamma, \varepsilon}(s ) 
	\bigr) 
	\Bigr )_{\! H_\Gamma } ds \qquad \label{p7}
\end{eqnarray} 
for all $t \in [0,T]$. We note that 
$$
	\int_0^t\!\!\! \int_{\Gamma}^{} \beta '_{\Gamma, \varepsilon }
	\bigl( u_{\Gamma, \varepsilon} (s) \bigr)
	\bigl | \nabla_\Gamma  u_\varepsilon (s) \bigr |^2 d\Gamma ds 
	\ge 0, 
	\quad 
	\varepsilon \int_0^t\!\!\! \int_{\Gamma }^{} u_{\Gamma, \varepsilon} (s) 
	\beta_{\Gamma, \varepsilon }\bigl( u_{\Gamma, \varepsilon} (s) \bigr) d\Gamma ds \ge 0
$$
{\pier due to the properties of $\beta_{\Gamma, \varepsilon }$.}
Then, {\pier recalling that 
$$\int_{\Gamma}^{} \widehat{\beta}_{\Gamma, \varepsilon}
	\bigl( u_{0\Gamma} \bigr) d\Gamma \leq 
	\int_{\Gamma}^{} \widehat{\beta}_{\Gamma}
	\bigl( u_{0\Gamma} \bigr) d\Gamma < +\infty
$$
by virtue of \eqref{p4}, and applying the Young inequality in the last term of \eqref{p7},}
we see that there {\pier exists} a positive constant $\tilde{M}_5${\pier ,
depending only} on $M_1$, $M_2$, $M_3$, $M_4$, $L_\Gamma$, $|u_{0\Gamma}|_{H_\Gamma }$ and $T$, {\pier such that} 
$$
	\bigl |\beta_{\Gamma, \varepsilon} (u_{\Gamma,\varepsilon }) 
	\bigr |_{L^2(0,T;H_\Gamma )} \le \tilde{M}_5.
$$
{\pier Hence, by comparison in \eqref{(30)} we also infer that}
$$
	|\Delta_\Gamma u_{\Gamma, \varepsilon} |_{L^2(0,T;H_\Gamma )} 
	 \le \tilde{M}_5
$$
and consequently (see, e.g.,\ \cite[Section~4.2]{Gr}])
$$
	|u_{\Gamma, \varepsilon }|_{L^2(0,T;H^2(\Gamma ))} 
	 \leq \left( |u_{\Gamma, \varepsilon }|_{L^2(0,T;V_\Gamma)}^2+
	|\Delta_\Gamma u_{\Gamma, \varepsilon }|_{L^2(0,T;H_\Gamma)}^2 \right)^{1/2}
	\le \bigl( M_1^2 T + \tilde{M}_5^2 \bigr) ^{1/2}
	 =: M_6. 
$$
{\pier In view of Lemma~3.4, using the theory of the elliptic regularity 
(see, e.g.{\pier ,\ \cite[Thm.~3.2, p.~1.79]{BrGi}}) along with 
the estimate $|u_{\Gamma, \varepsilon }|_{L^2(0,T;H^{3/2}(\Gamma ))}\le M_6$, 
it turns out that
$$
	|u_\varepsilon |_{L^2(0,T;H^2(\Omega ))} \le M_7
$$
for some constant $M_7$ independent of $\varepsilon \in (0,1]$.}\hfill $\Box$

\subsection{Passage to the limit}

In this subsection, we conclude the existence proof {\pier by passing to the limit on the 
sequence of approximate solutions. Indeed, owing to the estimates stated in the Lemmas
from~3.2 to~3.5, there exist} a subsequence of $\varepsilon $ {\pier (not relabeled)} and some limit functions $u$, $u_\Gamma$, $\xi$, $\xi_\Gamma$, {\pier $\lambda $} such that 
\begin{gather} 
	u_\varepsilon \to u \quad \mbox{weakly {\pier star} in } 
	H^1(0,T;H) \cap L^\infty (0,T;V) \cap L^2 \bigl (0,T;H^2(\Omega ) \bigr),  \label{(36)}
\\
	u_{\Gamma, \varepsilon} \to u_\Gamma \quad \mbox{weakly {\pier star} in } 
	H^1(0,T;H_\Gamma ) \cap L^\infty (0,T;V_\Gamma ) \cap L^2 \bigl (0,T;H^2(\Gamma ) \bigr),
	 \label{(37)}
\\ 
	\beta_\varepsilon (u_\varepsilon) \to \xi \quad \mbox{weakly in } 
	L^2(0,T;H),  \label{(38)}
\\ 
	\beta_{\Gamma ,\varepsilon } (u_{\Gamma, \varepsilon}) \to \xi_\Gamma 
	\quad \mbox{weakly in } 
	L^2(0,T;H_\Gamma ), \label{(39)}
\\
	\lambda _\varepsilon \to \lambda  \quad \mbox{weakly in } 
	L^2(0,T) \quad \label{(41)}
\end{gather} 
{\pier as  $\varepsilon \to 0.$ From \eqref{(36)} and \eqref{(37)}, due to 
strong compactness results  (see, e.g.,\ \cite[Sect.~8, Cor.~4]{Si})
we infer that 
\begin{gather} 
	u_\varepsilon \to u \quad \mbox{strongly in } 
	C([0,T];H) \cap L^2 (0,T;V), \label{(42)}
\\ 
	u_{\Gamma, \varepsilon} \to u_\Gamma  \quad \mbox{strongly in } 
	C([0,T];H_\Gamma ) \cap L^2 (0,T;V_\Gamma ). \label{(43)}
\end{gather} 
Moreover, on account of \eqref{(33)} it is a standard matter to deduce that
\begin{equation} 
	k_\varepsilon \to k \quad \mbox{weakly in } 
	H^1(0,T) \, \mbox{ and strongly in } C([0,T]),
 \label{(40)}
\end{equation} 
where 
$$ k_* \leq k (t):=\int_{\Omega }^{}w u (t) dx + \int_{\Gamma }^{} w_{\Gamma } u_{\Gamma}(t) d\Gamma  \leq k^* \quad \hbox{for all } \, t\in [0,T]. $$
We} point out that \eqref{(31)}, \eqref{(36)} and \eqref{(37)} imply that 
$u_\Gamma =u_{|_\Gamma}   $ a.e.\ on $\Sigma$, while \eqref{(32)}, \eqref{(42)} and \eqref{(43)} {\pier entail} 
$$
	u(0)=u_0 \quad \mbox{a.e.\ in } \Omega, 
	\quad 
	u_\Gamma (0) = u_{0\Gamma } 
	\quad \mbox{a.e.\ on } \Gamma.
$$
Now, \eqref{(41)}, \eqref{(40)} 
and the maximal monotonicity of $\partial I_{[k_*,k^*]}$ {\colli allow} us to conclude that 
$$
	\lambda \in \partial I_{[k_*,k^*]}(k)
	\quad \mbox{a.e.\ in } (0,T),
$$
{\pier while} \eqref{(42)}, \eqref{(43)} and 
the Lipschitz continuity of $\pi $ and $\pi_\Gamma $ imply that 
\begin{gather*} 
	\pi (u_\varepsilon) \to \pi (u) \quad \mbox{{\pier strongly} in } 
	C([0,T];H),
\\ 
	\pi_\Gamma (u_{\Gamma, \varepsilon}) \to \pi_\Gamma (u_\Gamma)  \quad \mbox{{\pier strongly} in } 
	C([0,T];H_\Gamma ) 
\end{gather*} 
{\pier as  $\varepsilon \to 0.$} 
At this point, we can pass to the limit in \eqref{(29)} and \eqref{(30)} obtaining 
\begin{gather*}
	\frac{\partial u}{\partial t}-\Delta u + \xi + \pi(u) + \lambda w = f
	\quad \mbox{a.e.\ in } Q,
\\
	\partial_\nu u
	+ \frac{\partial u_{\Gamma }}{\partial t}-\Delta_{\Gamma } u_{\Gamma }
	+ \xi_{\Gamma } +\pi_{\Gamma }(u_{\Gamma }) + \lambda w _{\Gamma } = f_\Gamma 
	\quad \mbox{a.e.\ on } \Sigma.
\end{gather*}
Let us comment that 
$\partial_\nu u_\varepsilon \to \partial_\nu u $ weakly in  
$L^2 {\fukao (}0,T;H^{1/2}(\Gamma ){\fukao )} $
as $\varepsilon \to 0 $,
due to \eqref{(36)} and the {\pier linearity and continuity} of 
the trace operator $u \mapsto \partial_\nu u$. 
Moreover, by applying \cite[p.~42, Proposition~1.1]{Ba} and using 
\eqref{(38)}--\eqref{(39)} with \eqref{(42)}--\eqref{(43)}, we obtain 
$$
	\xi \in \beta (u)
	\quad \mbox{a.e.\ in } Q,
\quad  \ 	\xi_\Gamma \in \beta_{\Gamma } (u_{\Gamma })
	\quad \mbox{a.e.\ on } \Sigma.
$$
Thus, it turns out that the pair 
$\mbox{\boldmath $ u$}=(u,u_\Gamma )$ is a 
solution of the limit problem, which can be stated exactly 
as in \eqref{(8)}--\eqref{13bis}. Also, we 
note the regularities 
$u \in C([0,T];V)$ and $u_\Gamma \in C([0,T];V_\Gamma )$ for the 
solution as a consequence of \eqref{(36)} and \eqref{(37)}. Morever, 
$\mbox{\boldmath $ u$}=(u,u_\Gamma )$ solves the abstract problem:
\begin{gather}
	\mbox{\boldmath $ u $} \in H^1(0,T;\mbox{\boldmath $ H $}) 
	\cap C([0,T];\mbox{\boldmath $ V $}), \label{p8}
\\
	\mbox{\boldmath $ u $}^* = (-\Delta u + \xi,\partial_\nu u -\Delta_\Gamma u_\Gamma +\xi_\Gamma  )
	\in L^2(0,T;\mbox{\boldmath $ H $}), \label{p9}
\\
	\lambda \in L^2(0,T),  \label{p10}
\\
	\mbox{\boldmath $ u $}'(t) + \mbox{\boldmath $ u $}^* (t) 
	+ \lambda (t) \mbox{\boldmath $ w $} + 
	\mbox{\boldmath $ \pi  $}
	\bigl (\mbox{\boldmath $ u $}(t) \bigr)  = \mbox{\boldmath $ f $}(t)
	\quad \mbox{in } \mbox{\boldmath $ H $}, 
	\ \mbox{for a.a.\ } t \in (0,T),\label{p11}
\\
	\mbox{\boldmath $ u $}^*(t) \in 
	\partial \varphi \bigl (\mbox{\boldmath $ u$}(t) \bigr ) 
	\quad \mbox{in } \mbox{\boldmath $ H$},
	\ \mbox{for a.a.\ } t \in (0,T),  \label{p12}
\\
	\lambda (t) \mbox{\boldmath $ w $} \in 
	\partial I_{\mbox{\boldmath \scriptsize $ K$}_{\mbox{\boldmath \tiny $ H $}}} \bigl (\mbox{\boldmath $ u$}(t) \bigr ) 
	\quad \mbox{in } \mbox{\boldmath $ H$}, 
	\ \mbox{for a.a.\ } t \in (0,T), \label{p13}
\\	\mbox{\boldmath $ u$}(0) = \mbox{\boldmath $ u$}_0 \quad \mbox{in } \mbox{\boldmath $ H$}.  \label{p14}
\end{gather}

\paragraph{Remark 3.3.} Let us point out that 
$$
	\mbox{\boldmath $ u$}^*(t)+\lambda (t) \mbox{\boldmath $ w$} \in \partial 
	( \varphi +I_{\mbox{\boldmath \scriptsize $ K$}})\bigl (\mbox{\boldmath $ u$}(t) \bigr )
	\quad \mbox{in } \mbox{\boldmath $ H$}, \ 
	\mbox{for a.a.\ } t \in (0,T).
$$
Therefore, {\colli \eqref{p8}--\eqref{p14} imply} that $\mbox{\boldmath $ u$}$ {\pier is the solution of the Cauchy problem expressed by} the abstract equation 
$$
	\mbox{\boldmath $ u $}'(t) 
	+ \partial (\varphi +I_{\mbox{\boldmath \scriptsize $ K$}})
	\bigl( \mbox{\boldmath $ u$}(t) \bigr)  + 
	\mbox{\boldmath $ \pi  $}
	\bigl (\mbox{\boldmath $ u $}(t) \bigr) \ni \mbox{\boldmath $ f $}(t)
	\quad \mbox{in } \mbox{\boldmath $ H $},
	\ \mbox{for a.a.\ } t \in (0,T)
$$
{\pier along with the initial condition \eqref{p14}.} 
Then, we emphasize that although the solution $\mbox{\boldmath $ u$}$ 
of {\pier this problem} is uniquely determined, the 
auxiliary quantities $\mbox{\boldmath $ u$}^*$ and $\lambda $ are 
not unique in general, {\pier except in special cases like, 
e.g.,} the case in which $\beta$ and  $\beta_\Gamma $ are single-valued 
(cf.\ \cite[Remark 2]{FuKe}).

\section{Continuous dependence}
\setcounter{equation}{0}

In this section, we prove Theorem 2.1.

\paragraph{Proof of Theorem 2.1.} 

Let {\pier $(\mbox{\boldmath $ u$}^{(1)}, \mbox{\boldmath $ \xi$}^{(1)}, \lambda^{(1)} )$
and 
$(\mbox{\boldmath $ u$}^{(2)}, \mbox{\boldmath $ \xi$}^{(2)}, \lambda^{(2)} )$
be two different solutions of {\rm (P)}, corresponding to  the} data 
($f^{(1)}$, $f^{(1)}_{\Gamma }$, $u_0^{(1)}$, $u_{0\Gamma }^{(1)}$) and 
($f^{(2)}$, $f^{(2)}_{\Gamma }$, $u_0^{(2)}$, $u_{0\Gamma }^{(2)}$), respectively. 
We take the difference between \eqref{p11} written for  $\mbox{\boldmath $ u$}^{(1)}(s)=(u^{(1)}(s),u_\Gamma ^{(1)}(s))$ 
and \eqref{p11} written for  $\mbox{\boldmath $ u $}^{(2)}(s)=(u^{(2)}(s),u_\Gamma ^{(2)}(s))$ at the time $t=s$ {\pier (note that the abstract equation \eqref{p11} 
comprehends both \eqref{(8)} and \eqref{(10)}). Then, we take 
the inner product with}
$\mbox{\boldmath $ u$}^{(1)}(s)-\mbox{\boldmath $ u$}^{(2)}(s)$ in 
$\mbox{\boldmath $ H$}$. {\pier Using the} monotonicity of 
$\beta $, $\beta_\Gamma$ and {\takeshi $\partial I_{\mbox{\boldmath \scriptsize $ K$}_{\mbox{\tiny \boldmath $ H$}}}$}, we 
obtain  
\begin{eqnarray*} 
	\lefteqn{ 
	\frac{1}{2} \frac{d}{ds} \bigl| u^{(1)}(s)-  u^{(2)}(s) \bigr|_H^2 
	+
	\frac{1}{2} \frac{d}{ds} \bigl| u^{(1)}_{\Gamma }(s)-  u^{(2)}_{\Gamma }(s) \bigr|_H^2 
	} \nonumber \\
	& ~ & {} +
	\bigl| \nabla u^{(1)}(s ) -\nabla u^{(2)}(s ) \bigr|_{H^{{\takeshi d}}}^2 
	+ 
	\bigl| \nabla _\Gamma u^{(1)}_{\Gamma }(s ) -\nabla _\Gamma u^{(2)}_{\Gamma }(s ) 
	\bigr|_{H_{\Gamma }^{{\takeshi d}}}^2 \\
	& \le & 
	\bigl( f^{(1)}(s ) -f^{(2)}(s ), u^{(1)}(s)-  u^{(2)}(s) \bigr)_{\!H} 
	+  
	\bigl( f^{(1)}_{\Gamma }(s ) -f^{(2)}_{\Gamma }(s ), 
	u^{(1)}_\Gamma (s)-  u^{(2)}_\Gamma (s) \bigr)_{\!H_\Gamma } \\
	& ~ & {}
	- \Bigl( \pi \bigl (u^{(1)}(s) \bigr) -  \pi \bigl (u^{(2)}(s) \bigr) , 
	u^{(1)}(s)-  u^{(2)}(s) \Bigr)_{\!H} \\
	& ~ & {}
	-  
	\Bigl( \pi_\Gamma  \bigl (u^{(1)}_\Gamma (s) \bigr) -  \pi_\Gamma  \bigl (u^{(2)}_\Gamma (s) \bigr), 
	u^{(1)}_\Gamma (s)-  u^{(2)}_\Gamma (s) {\fukao \Bigr)_{\! H}} \\
	& \le & \bigl| u^{(1)}(s)-  u^{(2)}(s) \bigr|_H^2 
	+ \bigl| u^{(1)}_\Gamma (s)-  u^{(2)}_\Gamma (s) \bigr|_{H_\Gamma }^2
	+  
	\frac{1}{2} \bigl| f^{(1)}(s ) -f^{(2)}(s )\bigr|_H^2 \\
	& ~ & {}+ \frac{1}{2}
	\bigl| f^{(1)}_{\Gamma }(s ) -f^{(2)}_{\Gamma }(s ) \bigr|_{H_\Gamma }^2 
	+ \frac{L^2}{2} \bigl| u^{(1)}(s)-  u^{(2)}(s) \bigr|_H^2 
	+ \frac{L_\Gamma ^2}{2}\bigl| u^{(1)}_\Gamma (s)-  u^{(2)}_\Gamma (s) \bigr|_{H_\Gamma }^2,
\end{eqnarray*} 
for all $t \in [0,T]$. {\pier Then, by integrating with respect to $s$ and applying the Gronwall lemma, it is straightforward to find a constant $C>0$, 
depending only on $L$, $L_{\Gamma}$ and $T$, such that 
the continuous dependence estimate \eqref{p0} holds.} \hfill $\Box$

{\pier
\section{Appendix}
\setcounter{equation}{0}

Let $\Omega \subset \mathbb{R}^{\takeshi d}$, ${\takeshi d}\geq 1$, be a bounded domain 
with smooth boundary $\Gamma:=\partial \Omega$. 
We use the same notation as in Section~\ref{main} for  $ \mbox{\boldmath $ H $} $ and   
$ \mbox{\boldmath $ V $} $.}	

\paragraph{Proposition 5.1.} {\it $\mbox{\boldmath $ V$}$ is dense 
in $\mbox{\boldmath $ H$}$.}

\paragraph{Proof.} For {\pier a fixed} $\mbox{\boldmath $ u$}=(u,u_\Gamma ) \in \mbox{\boldmath $ H$}$ and {\pier for} $n \in \mathbb{N}$, consider  the following {\pier elliptic} problem:
\begin{gather}
	v_n-\frac{1}{n}\Delta v_n =u 
	\quad \mbox{a.e.\ in } \Omega, 
	\label{(a1)}
\\
	\frac{1}{n} \partial_\nu v_n + (v_n)_{|_\Gamma}    =u_\Gamma 
	\quad \mbox{a.e.\ on } \Gamma.
	\label{(a2)}
\end{gather}
Then, {\pier let us write a variational formulation of \eqref{(a1)}--\eqref{(a2)}
\begin{align}
	\int_{\Omega }^{}v_n \eta \,dx + 
	\frac{1}{n} \int_{\Omega }^{} 
	\nabla v_n \cdot \nabla \eta \, dx 
	+ \int_{\Gamma }^{} (v_n)_{|_\Gamma}    \eta_{|_\Gamma}   d\Gamma 
	\nonumber \\ 
	= \int_{\Omega }^{} u \, \eta\, dx
	+ \int_{\Gamma }^{} u_{\Gamma}\,\eta_{|_\Gamma}   d\Gamma \quad \hbox{ for all } \eta    
	\in 
	V . \label{p15}
\end{align}
By applying the Lax-Milgram lemma, it is not difficult to see that for any 
$n \in \mathbb{N}$ there is a unique $ v_n \in  V $ that solves \eqref{p15}, i.e., 
satisfies} the above problem \eqref{(a1)}--\eqref{(a2)} with 
$\Delta v_n \in L^2(\Omega )$ and $\partial_\nu v_n \in L^2(\Gamma)$.
From the elliptic regularity for {\pier a} Neumann boundary condition 
(see, e.g.{\pier ,\ \cite[Thm.~3.2, p.~1.79]{BrGi}}), 
we {\pier infer} that $v_n \in H^{3/2}(\Omega )$, and this implies 
{\colli $(v_n)_{\Gamma} : = (v_n)_{|_\Gamma}    \in H^1(\Gamma)$.} 
Then, we have
$\mbox{\boldmath $ v$}_n = (v_n ,(v_n)_{\Gamma} ) \in \mbox{\boldmath $ V$}$
for all $n \in \mathbb{N}$.

{\pier Now, we take  $\eta = v_n$  in \eqref{p15} and apply the elementary 
Young inequality to deduce~that}
\begin{align}
	\frac{1}{2} \int_{\Omega }^{}|v_n|^2dx + 
	\frac{1}{n} \int_{\Omega }^{} 
	|\nabla v_n|^2dx 
	+
	\frac{1}{2} \int_{\Gamma }^{} \bigl| (v_n)_{\Gamma}    \bigr|^2 d\Gamma \nonumber
	\\
	\le \frac{1}{2} \int_{\Omega }^{} |u|^2 dx
	+ \frac{1}{2} \int_{\Gamma }^{} |u_\Gamma |^2 d\Gamma 
	=: M. \label{(a3)}
\end{align}
{\pier Hence, it turns out that
$\{v_n\}_{n \in \mathbb{N}}$ is bounded in $L^2(\Omega )$ and 
$\{ (1/\sqrt{n}\, ) \nabla v_n\}_{n \in \mathbb{N}}$  is bounded in 
$L^2(\Omega )^{\takeshi d}$,} with 
$$
	\left |\frac{1}{n} \nabla v_n \right |_{L^2(\Omega )^{{\takeshi d}}} \le \sqrt{\frac{M}{n}}
	\quad \mbox{for all } n \in \mathbb{N}. 
$$
Then, there exist {\pier a} subsequence $\{v_{n_k}\}_{k \in \mathbb{N}}$ 
of $\{v_n\}_{n \in \mathbb{N}}$ and ${\pier v} \in L^2(\Omega )$ such that 
\begin{equation}
	v_{n_k} \to {\pier v} \quad \mbox{weakly in } L^2(\Omega ),
\quad\ 
	\frac{1}{n_k} \nabla v_{n_k} \to 0 \quad \mbox{{\pier strongly} in } L^2(\Omega )^N \quad \mbox{as } k \to +\infty. 
	\label{(a4)}
\end{equation}
Next, choosing $\eta \in H^1_0(\Omega )$ in \eqref{p15}, we obtain
\begin{eqnarray*} 
	\int_{\Omega }^{}(v_{n_k} -u) \eta \, dx
     = - \int_{\Omega }^{} \frac{1}{n_k} \nabla v_{n_k} \cdot \nabla \eta \, dx 
	 \to 0 \quad \mbox{as } k \to +\infty,
\end{eqnarray*} 
namely, $v_{n_k} \to u$ in $H^{-1}(\Omega )$ as $k \to +\infty$; {\pier this means that the weak limit $v$ in \eqref{(a4)} should coincide with $u$ 
and the entire sequence 
\begin{equation}
\hbox{$v_{n}$ converges to $u $ weakly in $L^2(\Omega)$ as $n\to +\infty$.}
\label{(a5)} 
\end{equation} 
Now, from \eqref{(a3)} it follows that $\{ (v_n)_{\Gamma}    \}_{n \in \mathbb{N}}$ is bounded in $L^2(\Gamma )$; on the other hand, passing to the limit in \eqref{p15} we realize that  
$$ \lim_{n\to+\infty }  \int_{\Gamma }^{} \bigl( (v_n)_{\Gamma} - u_\Gamma  \bigr)  \eta_{|_\Gamma}   d\Gamma =0 \quad \hbox{for all } \eta \in H^1(\Omega), $$
whence 
\begin{equation}
	(v_{n})_{\Gamma} \to u_\Gamma \quad \mbox{weakly in } L^2(\Gamma ) \, \mbox{ as } n \to +\infty.
	\label{(a7)}
\end{equation}
Moreover, \eqref{(a3)} implies 
\begin{equation}
	\limsup_{n \to +\infty } \left\{ 
	\int_{\Omega }^{}|v_n|^2dx 
	+
	\int_{\Gamma }^{} \bigl| (v_n)_{\Gamma}\bigr|^2 d\Gamma 
	\right\} 
	\le \int_{\Omega }^{} |u|^2 dx
	+ \int_{\Gamma }^{} |u_\Gamma |^2 d\Gamma,
	\label{(a8)}
\end{equation} 
that entails the convergence of the norms 
$|v_n|_H$ and $| (v_n)_{\Gamma}|_{H_\Gamma}$
to $|u|_H$ and $| u_{\Gamma}|_{H_\Gamma}$, respectively. Thus, 
\eqref{(a5)}, \eqref{(a7)} and \eqref{(a8)} enable us to conclude that
$$
	v_{n} \to u
	\quad \mbox{strongly in } H, \quad 
	(v_{n})_{\Gamma} \to u_\Gamma
	\quad \mbox{strongly in } H_\Gamma,
$$
that is,
$$
	\mbox{\boldmath $ v$}_{n}=\bigl (v_{n},(v_{n})_{\Gamma}    \bigr ) 
	\to \mbox{\boldmath $ u$}=(u,u_\Gamma) 
	\quad \mbox{strongly in } \mbox{\boldmath $ H$} \, \mbox{ as } n \to +\infty, 
$$
which completes the proof.}\hfill $\Box$

\section*{Acknowledgments}

{\pier The} authors wish to express their heartfelt gratitude to 
{\pier Professors Goro Akagi and Ulisse Stefanelli}, 
who kindly gave them {\pier the opportunity} of exchange {\pier visits}
under the financial {\pier support of the JSPS-CNR joint research project} ``Innovative variational methods for evolution partial differential equations'', 2012-2013. 
{\pier The present note also benefits from a partial support of the MIUR-PRIN Grant 2010A2TFX2 ``Calculus of variations'' and the GNAMPA (Gruppo Nazionale per l'Analisi Matematica, la Probabilit\`a e le loro Applicazioni) of INdAM (Istituto Nazionale di Alta Matematica) for PC.}

\end{document}